\documentclass[twoside,12pt,reqno]{amsart}
\usepackage{upref,amsxtra,amssymb,amscd}
\usepackage{varioref}
\usepackage{verbatim}
\usepackage{epsfig}
\usepackage[french,english]{babel}
\usepackage{epic,eepic,eucal}
\usepackage[T1]{fontenc}
\usepackage{enumerate}

\def\a{         \alpha}

\newcommand{\GG}{{\mathcal T}}

\newcommand{\NN}{{\mathbb N}}
\newcommand{\RR}{{\mathbb R}}
\newcommand{\TT}{{\mathbb T}}
\newcommand{\ZZ}{{\mathbb Z}}

\newcommand{\QQ}{{\mathbb Q}}

\def\carre{ \hfill $\Box$    }

\newtheorem{theo}{\sc Theorem}[section]
\newtheorem{prop}[theo]{\sc Proposition}
\newtheorem{lemm}[theo]{\sc Lemma}
\newtheorem{coro}[theo]{\sc Corollary}

\theoremstyle{definition}

\newtheorem{defi}[theo]{\sc Definition}

\theoremstyle{remark}

\newtheorem{rema}[theo]{\sc Remark}

\numberwithin{equation}{section}

 \def\tn{|\kern-1pt|\kern-1pt|}

\begin{document}

\title{Rank one and mixing differentiable flows.}
\author{Bassam Fayad}
\maketitle

\begin{abstract} We construct, over some minimal translations of the two torus, special flows under a differentiable ceiling function that combine the properties of mixing and rank one. \end{abstract}


\section{\sc Introduction.} \label{vvvintroduction}

\subsection{Rank one and mixing.}  Rank one and mixing transformations or flows display the 
strong ergodic property of having minimal self joinings of all orders,
 a property which in  turn implies many features for the
transformation or flow such as having a trivial
centralizer and having  no factors\cite{Kingo,Kingttt,thouve,zeitz}.  Rank one and mixing transformations or flows are also mixing of any order \cite{Kalikow,Ryzhikov}.

\noindent The very few known examples of  transformations or flows combining the rank one property and mixing were all  produced in the same abstract frame of pure measure theory with  {\sl cutting and stacking} methods (see \S \ref{cut})
inspired by some works of Chacon and Ornstein \cite{Chacon,Ornstein} .
 While Chacon's seminal examples of cutting and stacking constructions were only weak mixing,  Ornstein  was the first to prove the existence of mixing rank one 
transformations. His existence result is based on probabilistic cutting  and stacking constructions with random
spacers (stochastic constructions of mixing rank one flows were obtained in \cite{Prikhodko}). Later, explicit cutting and stacking constructions were proven to be of rank one and mixing \cite{adams}.

In this paper we give a differentiable realization of rank one and mixing in the case of flows.  More precisely, we will define in Section  \ref{translation}
(following \cite{Y}) an uncountable
dense subset $Y \subset {\RR}^2$, for which we will prove the following 
\begin{theo} \label{theoreme} For any $(\a,\a') \in
Y$, there exists 
a strictly positive real function 
$\varphi$ defined on ${\TT}^2$ of class $C^1$  such that  the  special flow built over
$R_{\a,\a'}$ with the ceiling function $\varphi$ is of  rank one and mixing (with respect to its unique invariant measure). 
\end{theo} 

\noindent The exact definitions of special flows, mixing and rank one for transformations and for flows are given in Section \ref{preliminaries}. Roughly, a
 measure preserving transformation or flow on $(X, {\mathcal A}, \nu)$ is said to be of rank one if it has a sequence of towers that asymptotically generate the $\sigma$-algebra ${\mathcal A}$ (see Section \ref{r1}).   The
property of rank one  essentially reflects the existence of some cyclic
approximations for the flow; and clearly, cyclic
approximations do not favor mixing. For instance,  studying cyclic
 approximations as defined by  Katok and Stepin shows that the speed of approximation required to guarantee a
simple spectrum (a weaker property than rank one) implies the absence of mixing \cite{KS}.

On the other hand, we know, since the work of Katok, Stepin  and Shklover, that special flows over Liouvillean irrational rotations of the torus, even with  analytic ceiling functions, can combine fast cyclic
approximations (implying rank one) \cite{Ka1}  as well as some mixing features, namely weak mixing  \cite{katokrobinson,Kkk,S}. Using Baire category arguments, rank one and weak mixing can also be derived for most of the time-maps of these special flows \cite{weakmixing}. But these examples of weak mixing special flows and transformations are nevertheless {\sl rigid}, in the sense that $T^{t_n} \rightarrow Id_{{\TT}^2}$ for some sequence $t_n \rightarrow \infty$.  Rigidity of smooth special flows over irrational rotations of the circle is due to an improved Denjoy-Koksma inequality 
 involving the Birkhoff sums of the ceiling  function over the rotation.

\subsection{Uniform stretch mixing for special flows over translations.} To obtain mixing special flows over translations one can either consider special flows over circular rotations and under ceiling functions with singularities \cite{ko2,Khanin,bsmf} or turn  to some minimal translations on higher dimensional tori for which the Denjoy-Koksma inequality does not hold \cite{Y} and over which it is possible to construct mixing special flows with real analytic ceiling functions \cite{mixing}. In both cases, the key underlying mixing is the uniform stretch at all times of the Birkhoff sums of the ceiling function (see Section \ref{uuss} below). Under this stretch,  the image of a small interval on the base  becomes  as time goes to infinity increasingly  and uniformly expanded along the fibers
of the special flow hence tending to be 
equally distributed in the space by unique ergodicity of the translation on the base.  However, all these examples  most likely fail to be of finite
rank. Indeed, as we will observe later, the {\sl uniform stretch property at all times plays against the rank one  property.}

\subsection{The cutting and stacking techniques.}  \label{cut}

In the construction of rank one transformations by the cutting and stacking techniques the space and the transformation are obtained in the same time by considering successive towers  of intervals as towers for the transformation (see Section \ref{r1}). Each tower $C_{n+1}$ is obtained by cutting the tower $C_n$ into $r_n$ subcolumns  of equal width and adding some number $l_{n,i}$ of {\sl spacers} above every $i^{th}$ subcolumn before stacking over it the $(i+1)$st subcolumn for $i=1,...,r_n-1$.

 With the latter cutting and stacking construction, it is possible to avoid the cyclic approximations 
and obtain mixing transformations by actually blowing
up the top levels of the successive rank one towers.  But in this case,  as long as the levels of the towers are thought of as 
intervals, it appears difficult to adapt the constructions into a differentiable frame.

One of the crucial points in  our construction is the following
elementary fact: {\sl small measure in dimension $2$ is not equivalent
to small diameter}. More precisely, assume we are stacking in a column  
 $n^2$ disjoint squares of equal area  ${1 / n^2}$, doing so in an
 isometric way until we reach the top  $n^{3 \over 2}$ squares, and
 then applying  uniform stretch of magnitude 
${1 / n^{4 \over 3}}$ on each of the last $n^{3 \over 2}$ levels
(clearly this can be done with small derivatives since $1 / n^{4 \over 3} =o ( 1 / n) $); we
will thus reach the top of the tower  with a cumulated stretch
amounting to   
$n^{{3 \over 2}-{4 \over 3}}$ which is large  while the total measure of the
squares where the stretch was applied is   ${n^{{3 \over 2} -2}} = o(1)$.

\subsection{Combining uniform stretch and cutting and stacking.}

The construction we will present here of a  mixing rank one flow combines the uniform stretch and the cutting and stacking techniques.

Let us first give a brief description of the mixing special flows constructed in \cite{mixing}: They are special flows above a minimal translation $R_{\a,\a'}$ of $\TT^2$ and under a ceiling function given by
$$\varphi(x,y) = 1+ \sum_{n =2}^{\infty} {\cos(2 \pi q_n x) \over e^{q_n}} + {\cos(2\pi q'_ny) \over e^{q'_n}}$$
where ${\lbrace q_n \rbrace }_{n \in \NN}$ and ${\lbrace q'_n \rbrace }_{n \in \NN}$ are the sequences of denominators of the convergents of $\a$ and $\a'$. If these sequences are such that $q_n' \geq e^{3q_n}$ and $q_{n+1} \geq e^{3q'_n}$ for all $n \in \NN$, then due to the term $\cos(2 \pi q_n x) / e^{q_n}$ we obtain that $S_m \varphi (x,y)$ is uniformly stretching in the $x$ direction for $m \in [e^{2q_n}, e^{2q'_n}]$ while the term 
 $\cos(2 \pi q_n' y) / e^{q'_n}$ is responsible for uniform stretch in the $y$ direction of $S_m \varphi (x,y)$ for $m \in [e^{2q'_n}, e^{2q_{n+1}}]$. Since the latter intervals cover  a neighborhood of infinity  $[n_0, \infty) \subset \NN$ we deduce that the special flow is mixing.

Here, we will use essentially the same translations on the base but $\varphi$ must be modified in order
to gain the property of rank one for the special flow without losing the mixing property.
 The modification is done in the $C^1$ topology but can be made smoother on higher dimensional tori.

First of all, a criterion that guarantees the  rank one property for a special flow over a rank one transformation is given: Starting with a  rank one sequence of
towers of the transformation on the base (see Section \ref{r1}), 
the idea of the
criterion is that under a condition of flatness of the Birkhoff sums of the ceiling function computed over the base  of the successive towers, it is possible to {\sl lift} these  towers into rank one towers for the special flow. Hence, we first choose for the translation on the base  a particular sequence of rank one towers that we want to lift into rank one towers for the special flow and construct the ceiling function subsequently. 

Next, starting with the function $\varphi$ as above, when  the uniform stretch of a term like $\cos (2 \pi q_n x) / e^{q_n}$  or  $\cos (2 \pi q'_n x) / e^{q'_n}$ jeopardizes the rank one property, we have to change it. We do this as follows:

\begin{itemize} 

\item We can change the cosine by a function that is essentially flat everywhere except over the top steps of the rank one towers that we want to lift, confining thus the uniform stretch to the top levels of the towers. This can be done smoothly as was explained above. 

\item We can replace the cosine by a staircase-like function constant on the levels of the rank one towers and smoothened-up with bump functions. 

\end{itemize}

 By the first procedure, rigidity times due to fast cyclic approximations are precluded and replaced by mixing sequences of time due to the uniform stretch. However before reaching the top levels of a tower of a translation there are intermediate rigidity times that will also be rigidity times for the special flow due to the flatness of the Birkhoff sums.

 By the second procedure, uniform stretch of the Birkhoff sums that provided mixing via  uniform continuous stretch  of intervals  gives way to non-uniform stretch, i.e. staircase stretch of intervals similar to the one obtained on the top of each column in the cutting and stacking constructions. Subsequently, the proof of mixing involves arithmetically spaced Birkhoff averages under the action of the flow. These averages are shown to converge using the mixing times obtained previously from uniform stretch.

\subsection{Plan of the construction.} In the next section we introduce some definitions and notations and we state a criterion that guarantees the  rank one property for a special flow over a rank one transformation. We also recall the criterion that yields mixing  for the special flow from uniform stretching of the Birkhoff sums of the ceiling function. 

Section \ref{translation} is reserved to the choice of the translation $R_{\a,\a'}$ and to the description of a special rank one sequence of  towers for $R_{\a,\a'}$ that we will want to lift into a rank one sequence  of  towers for the flow.

In Sections \ref{pX} and \ref{pY}, we list the properties  required on the functions $X_n(x,y)$ and $Y_n(x,y)$ that will substitute the terms $\cos (2 \pi q_n x) / e^{q_n} $ and $\cos (2 \pi q'_n y) / e^{q'_n}$ in the expression of the ceiling function. The effective construction of these functions is done in the last sections \ref{chiconstruct} and \ref{psiconstruct}. 

In Sections \ref{proof of rank one} and \ref{proof of mixing}  respectively, we prove  that the special flow built over $R_{\a,\a'}$ and under the  function $\varphi = \varphi_0 + \sum_{n \geq n_0} X_n + Y_n $ is of rank one and mixing, where $\varphi_0$ and $n_0$ are chosen so that $\varphi$ is strictly positive and of mean value $1$.

\section{Preliminaries} \label{preliminaries}

\subsection{Special flows.} \label{special flow} Let $(M, T, \mu)$ be a dynamical system. Given a real function $\varphi \in
L^1(M,\mu)$ such that $\varphi \geq c > 0$, we define the {\sl special flow constructed over $T$ and under the ceiling function $\varphi$}  to be
the quotient flow of the action $M \times \RR  \rightarrow  M \times
\RR :    (z,s)  \rightarrow  (z,s+t) $
by the relation $(z, s + \varphi(z)) \sim (Tz, s)$. This
flow acts on the manifold $M_{T,\varphi}$
obtained from  the subset of $M \times \RR$: $M_{\varphi} = \left\{ (z,s) \in  M \times \RR \ / \ 0 \leq s
< \varphi ({z}) \right\}$ by identifying pairs $(z,\varphi(z))$ and
 $(Tz,0)$. It  
 preserves the normalized product measure on $M_{T, \varphi}$, i.e. the
product of the measure $\mu$ on the base  with the Lebesgue
measure on the fibers divided by  $\int_{M}
\varphi(z) d\mu(z)$.
Moreover, if the transformation $T$ is uniquely ergodic then so is the special flow. 
We denote by $T^t_{T, \varphi}$ the special flow above a transformation $T$ and under the ceiling function $\varphi$. The Birkhoff sums of $\varphi$ over the iterates by $T$ of a point $z \in M$ are denoted by $ S_m \varphi(z) = \sum_{i=0}^{m-1} \varphi ( T^i z ). $

\begin{defi} \label{mzt} For $(z,s) \in M_{T, \varphi}$ and $t \in \RR^+$ we introduce the notation $m(z,s,t)$ for the only $m \in \NN$ that satisfies 
\begin{eqnarray*}   0 \leq s + t - S_m \varphi(z) \leq \varphi(T^mz). \end{eqnarray*} 
With this definition of $m=m(z,s,t)$ we get 
$$T^t(z,s) = \left( T^mz, s+t- S_{m} \varphi(z) \right).$$
\end{defi}
For a point $z \in M$ we sometimes use the notation $z$ and $m(z,t)$ for $(z,0) \in M_{T, \varphi}$ and $m(z,0,t)$.

\subsection{The rank one property.} \label{r1}

\subsubsection{Rank one transformations.} Let $(M,T, \mu)$ be a dynamical system. Given a measurable set $A \subset M$ and an integer $h$ such that $A, T(A), ..., T^{h-1}A$ are disjoint we say that  $A \sqcup TA \sqcup ... \sqcup T^{h-1}A$  is a {\sl tower} of $T$ and denote it by $\GG (A, h) $. The set $A$ is called the {\sl base  } of the tower and $h$ its {\sl height}.  Every $T^{k}(A)$, $ k \leq h-1$, is called a {\sl level} of the tower. The {\sl measure} of the tower is $h \mu(A)$. Rokhlin lemma insures that if the set of periodic points of $T$ is of measure zero then for any $\epsilon > 0$ and any $h \in \NN^*$ there exists a tower of $T$  with height  $h$ and measure greater than $1 - \epsilon$.

 Given a measurable partition ${\mathcal P}$ of $M$ and an $\epsilon > 0$, we say that a measurable set $A$ is $\epsilon-${\sl monochromatic} with respect to ${\mathcal P}$ if all but a proportion less than $\epsilon$  of the set $A$ is included in one atom of ${\mathcal P}$. We say that a tower of $T$ with base $A$ and height $h$ is $\epsilon-${\sl monochromatic} with respect to ${\mathcal P}$ if all but a proportion less than $\epsilon$ of its levels (i.e. less than $h \epsilon$ levels) are $\epsilon-${\sl monochromatic} with respect to ${\mathcal P}$.
 
\begin{defi} We say that a dynamical system  $(M,T,\mu)$ is of {\sl rank one} or has the {\sl rank one property} if for any finite measurable partition ${\mathcal P}$ of $M$ and for any
$ \epsilon >0  $  there exists a tower for $T$ that has measure greater than $1 - \epsilon$ and is $\epsilon-${\sl monochromatic} with respect to ${\mathcal P}$.  In other words, $(M,T,\mu)$ is of rank one if there exists a sequence of towers  for $T$  that generates the sigma algebra of
finite partitions of  $(M, \mu)$. \end{defi}

\subsubsection{Rank one flows.} Let  $ (M,T^t,\mu)$ be an ergodic flow. For any positive real number $H$ and any $\epsilon>0$, we can represent  ${\lbrace T^t \rbrace}_{t \in \RR}$ as a special flow over a system $(X,T,\nu)$ with a ceiling function $\varphi$ such that: 
\begin{itemize}
\item[(i)] $\varphi(x) \leq H $ for every $x$ in $X$,
\item[(ii)] $\varphi(x) = H$ on a subset $B \subset X$ of measure $ \nu(B) \geq 1
- \epsilon.$ \end{itemize}
 This flow version of Rokhlin's lemma was first introduced by
Ornstein in \cite{orsos}. We call  $\bigcup_{t=0}^{H} T^tB  $ a { tower} of  
${\lbrace T^t \rbrace}_{t \in \RR}$
with base  $B$ and height $H$. Every $T^t(B)$, $t\leq H$, is a {\sl horizontal level} of the
tower. The measure of the tower is $H \nu(B) / \int_{X} \varphi(x) d \nu(x) \geq \nu(B) \geq 1- \epsilon$.

Given a finite measurable partition $\mathcal P$ of  $(M,\mu)$ 
and an $\epsilon > 0$, we say that a horizontal level $T^sB$, $s \leq H$
of the tower is $\epsilon$-monochromatic with respect to $\mathcal P$ if a proportion not less than $1-\epsilon$ of the $T^s_* \nu$-measure of
 this horizontal level is included in one of the atoms of $\mathcal P$.  We say that a tower  above $B$ of height $H$ is 
 $\epsilon$-monochromatic with respect to $\mathcal P$ if all but a proportion less than $\epsilon$ of its horizontal levels (proportion measured with respect to the Lebesgue measure on $[0,H]$) are 
$\epsilon$-monochromatic with respect to $\mathcal P$.

Like the Rokhlin lemma, the definition of  rank one can also be stated for flows as in \cite{zeitz}
\begin{defi} Let  $(M,T^t, \mu)$ be a dynamical system. We
say that the flow ${\lbrace T^t \rbrace}_{t \in \RR}$ is of rank one  if for any finite measurable partition $\mathcal P$ of $M$ and for every
$ \epsilon  >0$, there is a tower for ${\lbrace
T^t \rbrace}_{t \in \RR}$ of measure greater than $1 - \epsilon$  that is $\epsilon-$monochromatic with
respect to $\mathcal P$. 
\end{defi}

\subsubsection{A criterion that guarantees the property of rank one for special flows over rank one transformations.} 
The criterion involves the Birkhoff-sums of the ceiling function and allows us to "lift"  
 rank one towers of the transformation on the base to rank one towers  of the flow.  

Let $(M,T, \nu)$ be a dynamical system of rank one.  We call a sequence of towers of $T$, ${ \lbrace \GG ( B_n, h_n) \rbrace}_{n \in \NN}$, a {\sl rank one sequence of towers} for $T$ if given any $\epsilon>0$ and any finite measurable partition ${\mathcal P}$ of $M$, there exists $n_0$ such that for every $n \geq n_0$, 
  the tower ${\GG (B_n, h_n)}$ has measure greater than $1 - \epsilon$ and is $\epsilon$-monochromatic with respect to ${\mathcal P}$. 
 
\begin{prop}[{\bf Criterion for rank one}] \label{criterion}
 Let $(M,T, \nu)$ be a dynamical system of rank one and 
 $T^t_{T, \varphi}$ be a special flow  constructed over $T$ with a  ceiling function satisfying $0< c \leq \varphi \leq C < \infty$. If there exists a rank one sequence of towers,  ${\lbrace \GG (B_n, h_n) \rbrace}_{n \in \NN}$ of $T$ such that 
\begin{eqnarray} \label{rangl}
\sup_{m\leq h_n} \sup_{z,z' \in B_n} |S_m \varphi(z') - S_m \varphi(z)| \mathop{\longrightarrow} \limits_{n  \rightarrow \infty}  0,
\end{eqnarray} 
then the flow $T^t_{T, \varphi}$ is of rank one. 
\end{prop}

\noindent{\sl Proof.} Since $\inf_{z \in M} \varphi(z) \geq c >0$,  (\ref{rangl}) implies that for sufficiently large $n$ there exists $H_n \in \RR$ such that    $ \sup_{z \in B_n} S_{h_n-1} \varphi(z) < H_n <\inf_{z \in B_n} S_{h_n} \varphi (z)$. Hence in Definition \ref{mzt} we have for any $z \in B_n$ $m(z,H_n)= h_n -1$, therefore $\bigcup_{t=0}^{H_n} T^tB_n$  is a tower of the special flow since $ \GG (B_n, h_n)  $ is a tower of $T$. Moreover since  the measure of the tower on the base $\GG ( B_n, h_n -1 ) $ is greater than $1 - \epsilon -\nu (B_n) $ we get that the measure of the flow tower $\bigcup_{t=0}^{H_n} T^tB_n$ is greater than $1 - (\nu(B_n)+ \epsilon) \sup_{z \in M} \varphi(z) / \int_M \varphi(z) d \nu (z)$. 
 
On the other hand fix an arbitrary measurable set $A \subset M$. From (\ref{rangl}) we have that as the levels of the towers ${\GG ( B_n, h_n)}$
 of $T$ become increasingly monochromatic with respect to the set $A$ then the levels of the flow towers 
 $\bigcup_{t=0}^{H_n} T^tB_n$ become increasingly monochromatic with respect to  $\bigcup_{t=s_0}^{s_0+s_1}A$ where $s_0, s_1 \in \RR$. Since any finite measurable partition of $M_{T,\varphi}$ can be approximated by sets of the latter form the proof of rank one for the flow follows.    \carre

\subsection{The mixing property.} \label{mixing property}

\subsubsection{Mixing sequences of time.}  \label{mixing sequences} We  recall that a dynamical system $(M, T^t, \mu)$ is said to be {\sl mixing}  if for any measurable sets $A, B \subset M$ one has
\begin{eqnarray} \label{m1}
\lim_{t\rightarrow \infty} \mu ( T^{-t} A \cap B )  = \mu(A) \mu (B). 
\end{eqnarray}

\begin{defi} \label{m132} We say that a sequence $t_n \rightarrow \infty$ is a {\sl mixing sequence} for the flow ${\lbrace T^t \rbrace}_{t \in \RR}$ if (\ref{m1}) holds along the sequence $t_n$ as $n$ goes to infinity. A sequence of subsets of $\RR$, ${\lbrace I_n \rbrace }_{n \in \NN}$ is called a {\sl mixing sequence of sets} for the flow ${\lbrace T^t \rbrace}_{t \in \RR}$ if any sequence $t_n \in I_n$ is a mixing sequence. 
\end{defi}

Clearly, if a neighborhood of $\infty$  can be covered by a finite union of mixing sequences of sets, e.g. for some $a \in \RR$, $[a, \infty) \subset \bigcup_{j=1}^{j_0} \bigcup_{n \in \NN} I_{j,n}$, where each ${\lbrace I_{j,n} \rbrace}_{n \in \NN}$ is a  mixing sequence of sets then the flow is mixing.

\subsubsection{Good sequences of partial partitions.} \label{good}

 We denote by a partial partition of $(M, \mu)$ a finite  collection of  disjoint measurable sets in $M$. In the case $M=M_{R_{\a,\a'}, \varphi}$  we also call partial partitions collections of disjoint sets of $M$ of the form $ I \times \lbrace y,s \rbrace$, $(y,s) \in \TT \times \RR$, $I$ interval on $\TT$,  or collections of disjoints sets of the form $ R \times \lbrace s \rbrace $, $s \in \RR$, $R$ rectangle in $\TT^2$.

\begin{defi} Let  ${\lbrace {\Omega_t}  \rbrace}_{t \in \RR}$  (or $t \in \NN$) be a family of partial partitions of $(M, \mu)$. We say that  ${{\Omega} }_{t}$  tends to the partition into points as $t$  goes to infinity and write ${\Omega}_t \rightarrow \epsilon$ if every  measurable subset of $M$ becomes arbitrarily well approximated in measure by unions of sets in ${\Omega}_t$ as $t$ goes to infinity. \end{defi}

\begin{defi}  \label{fde} Let $(M,T,\mu)$ be a dynamical system. Let  ${\lbrace t_n \rbrace}_{n \in \NN}$ be a sequence of real numbers and $A$ a measurable subset of $M$. We say that a  family of partial partitions  ${\lbrace {\Omega}_n \rbrace }_{n \in \NN}$ is {\sl good for the sequence ${\lbrace t_n \rbrace}_{n \in \NN}$   and for  $A$} if for any $\varepsilon > 0$ there exists $n_0 \in \NN$ such that  for any integer $n \geq n_0$, for any atom $\xi \in {\Omega}_{n}$ 
\begin{eqnarray}   \left| \mu ( \xi  \cap T^{-t_n} A) - \mu (\xi)  \mu ( A ) \right| \leq  \varepsilon \mu (\xi). \label{fubfub} \end{eqnarray}
If for a sequence  ${\lbrace t_n \rbrace}_{n \in \NN}$, we have for any  measurable set $A$ a sequence $\Omega_n(A) \rightarrow \epsilon$  that is good then  ${\lbrace t_n \rbrace}_{n \in \NN}$ is a mixing sequence for the flow ${\lbrace T^t \rbrace}_{t \in \RR}$. 

In the case of partial collections $\Omega_n$ of $M=M_{R_{\a,\a'}, \varphi}$
 with positive codimension sets we will also say that $\Omega_n$ is good for  ${\lbrace t_n \rbrace}_{n \in \NN}$  
and for a measurable set $A \subset M$ if for any $\varepsilon > 0$ there exists $n_0 \in \NN$ such that  for any integer $n \geq n_0$, for any atom $\xi \in {\Omega}_n$
we have
\begin{eqnarray}   \left| \lambda^{(i)} ( \xi  \cap T^{-t_n} A) - \lambda^{(i)} (\xi)  \mu ( A ) \right| \leq  \varepsilon \lambda^{(i)}  (\xi), \label{fub2} \end{eqnarray}
where $\lambda^{(1)} = \lambda$ is the Lebesgue measure on the line and $\lambda^{(2)} = \lambda \times \lambda$ and $i =1,2$ depending on the dimension of $\xi$. A Fubini argument then gives the same conclusions related to mixing as above.
\end{defi}

\begin{rema} In the case of a metric space $M$ and a Borelian measure $\mu$, mixing for a sequence ${\lbrace t_n \rbrace }_{n \in \NN}$      follows if we check the conditions of Definition \ref{fde} for any ball $A \subset M$. \end{rema}

Along the line of the definitions above, we have that if for any measurable set $A$ there exists a family of partial partitions  ${\lbrace {\Omega}_t  \rbrace}_{t \in \RR}$  such that $\Omega_t \rightarrow \epsilon$ and that ${\lbrace {\Omega}_{t_n}  \rbrace}_{n \in \NN}$ is good for $A$ and ${ \lbrace t_n \rbrace }_{n \in \NN}$ as long as $t_n \rightarrow \infty$ then the system $(M,T^t, \mu)$ is mixing.

\subsubsection{Uniform stretch.} \label{uuss} One of the tools we will use to derive mixing is uniform stretch. We recall the definition for a real function on a segment  $[a,b] \subset \RR$ (see \cite{ko2,mixing})

\begin{defi} Let $\varepsilon > 0$ and $K>0$. We say that a real function $g$ defined on an interval $[a,b]$ is {\sl
$(\varepsilon, K)$-uniformly stretching} on  $[a,b]$ if 
$$ \sup_{[a,b]} g - \inf_{[a,b]} g \geq K,$$ 
and if for any  $\inf_{[a,b]} g \leq u \leq
v \leq \sup_{[a,b]} g $, the set
$$I_{u,v} = \lbrace x \in [a,b] \ / \ u \leq g(x) \leq v \rbrace, $$
has Lebesgue measure 
$$ (1- \varepsilon){v-u \over |g(b)-g(a) |} (b-a) \leq \lambda(I_{u,v}) \leq (1+ \varepsilon){v-u \over |g(b)-g(a)|}(b-a).$$
\end{defi}

We assume now that $g$ is at least two times differentiable  and we recall the following   straightforward but useful criterion on the derivatives
of $g$ insuring its uniform stretch on the segment $[a,b]$:   

\begin{lemm} [A  Criterion  for uniform stretch] \label{deriveeseconde}
If 
\begin{eqnarray*}
 \inf_{x \in [a,b] } | g'(x) | |b-a| &\geq& K \\
 {\rm and} \ \ \ \sup_{x \in
[a,b]}|g''(x)| |b-a| &\leq& \varepsilon \inf_{x \in [a,b]} |g'(x)|
\end{eqnarray*} 
then $g$ is   $(\varepsilon, K)$-uniformly
stretching on  $[a,b]$.
\end{lemm}

The following proposition derives  mixing from uniform stretch for a special flow above a minimal translation of the two torus $\TT^2 = \RR^2 / \ZZ^2$. Its proof can be found in \cite{mixing}.

\begin{prop} [{\bf Criterion for a mixing sequence}] \label{criterion mixing} Let $T^t_{R_{\a,\a'}, \varphi}$ be the special flow constructed over some minimal translation $R_{\a,\a'}$ of $\TT^2$ and under a strictly positive continuous  function  $\varphi$ two times differentiable in the $x$ direction. Let ${\lbrace t_n \rbrace }_{n \in \NN}$ be a sequence  of  real numbers. If  ${\lbrace {\Omega}_n \rbrace}_{n \in \NN}$   is a sequence of partial partitions of the circle in intervals  and if there exist sequences $K_n \rightarrow \infty$ and $\nu_n \rightarrow 0$ such that
\begin{itemize} \item For any interval $I \in {\Omega}_n$, for any $y \in \TT$, for $m=m(x,y,t_n)$ for some $x \in I$, $ S_m \varphi(.,y)$ is $(K_n, \nu_n)$-uniformly  stretching  on   $I \times \lbrace y \rbrace$, 
\end{itemize} 
then any sequence of partial partitions of $M_{R_{\a,\a'},\varphi}$ with sets of the form $\xi \times \lbrace y,s \rbrace$ , $\xi \in {\Omega}_n  $  is  good for ${\lbrace t_n \rbrace }_{n \in \NN}$ (and for any measurable set $A$). Therefore, if in addition ${\lbrace {\Omega}_n \rbrace}_{n \in \NN} \rightarrow \epsilon$ then $
 {\lbrace t_n \rbrace }_{n \in \NN}$ is a mixing sequence for the special flow.
\end{prop}

\subsection{Notations.} \label{notations}

\subsubsection{} For $d \in {\NN}^*$ let $\TT^d = \RR^d / \ZZ^d$. For $k \in \NN \cup \lbrace + \infty \rbrace $ we denote by
$C^{k} ({\TT}^d, \RR)$ the set of real functions on ${\RR}^d$ of class $C^{k}$
and ${\ZZ}^d$-periodic. By  $C^{k} ({\TT}^d, {\RR}^*_+)$ we will
denote the subset of $C^{k} ({\TT}^d, \RR)$ of strictly positive functions.
We will use the notation  $ \|\varphi \| := \sup_{z \in {\TT}^d} |\varphi(z)|. $
 For $r \in \NN$, the notation $D_x^r \varphi$ is used
for the derivative of order $r$ of $\varphi$ with  respect to $x$ and the norm   on $C^{k} ({\TT}^d, \RR)$ we consider is ${\|\varphi\|}_{C^k}:= \sum_{r+p=k} \|D_x^rD_y^p \varphi\|$.

\subsubsection{} Let $x$ be a real number; we denote by: $ $
\\--$[x]$ the  integer part of $x,$
\\--$\lbrace x \rbrace = x - [x]$ its fractional part, $ $
\\--$\tn x \tn = \min(  \lbrace x \rbrace, 1 - \lbrace x
\rbrace)$ the distance of $x$ to the closest integer. $ $

\subsubsection{} We recall some facts about the best approximations of an irrational number by rational ones. When we write ${p \over q} \in \QQ$ we assume that $q \in \NN$,
$q \geq 1$, $p \in \ZZ$ and that $p$ and $q$ are relatively prime. To each $\a \in \RR \setminus \QQ$, there exists a sequence of  rationals
${\lbrace {p_n /
  q_n} \rbrace}_{n \in \NN} $ called the {\sl convergents} of $\a$, such that:
\begin{eqnarray} \label{reduites} \tn q_{n-1} \a \tn < \tn q \a \tn  \ \ {\rm for \ every \ }  0 < q < q_{n}, \ q \neq q_{n-1} \end{eqnarray} and for any $n \in \NN$ 
\begin{eqnarray} \label{reduites1} { 1  \over q_n (q_n + q_{n+1})} \leq {(-1)}^n ( \a - {p_n \over q_n}) \leq   { 1  \over q_n q_{n+1}}. \end{eqnarray}

Let $\a' \in \RR \setminus \QQ$ and assume that $\a' - {p'_{{n}-1} \over q'_{{n}-1}} > 0 $  then we have 
\begin{eqnarray}
 \left| \a' - {p'_{{n}-1} \over q'_{{n}-1}} - { 1 \over q'_{{n}-1} q'_{n} } \right| < { 1
\over q'_{n} q'_{{n}+1}}. \label{reduites2} \end{eqnarray}
 
\begin{rema}  In all the paper we will always assume a fixed parity for $n$, say $n$ odd, so that $ \a' - {p'_{{n}-1} \over q'_{{n}-1}} > 0 $. The constructions we would have to make at step $n$ if $n$ is even being  similar to the ones we will make assuming $n$ is odd.   \end{rema}

\section{\sc The translation on the base ${\TT}^2$.} \label{translation}

\subsection{The choice of the translation $R_{\a,\a'}$.} \label{choice}

Following \cite{Y} we introduce
\begin{defi} \label{yocco} Let $Y$  be the set of couples $(\a,\a')
\in {{\RR}^2 \setminus \QQ}^2$ such that the sequences of denominators of the convergents of $\a$ and
$\a'$,  ${\lbrace q_n \rbrace}_{n \in \NN}$ and  ${\lbrace q'_n \rbrace}_{n \in \NN}$ respectively, satisfy the
following: there exists $n_0 \in \NN$ such that, for any $n \geq n_0$
\begin{eqnarray}
\label{10000} &q'_n \ \ \geq  e^{3q_n}, \\
\label{20000} &q_{n+1} \geq e^{3q'_n}, \\
\label{30000} &q_n \wedge q'_{n-1} = 1, \ {\rm and } \ q_n \wedge q'_n = 1. \end{eqnarray}
\end{defi} 
\noindent Here $p \wedge q =1 $ stands for $p$ and $q$ relatively prime.

 The importance
of (\ref{10000}) and (\ref{20000}) in the choice of $(\a,\a')$ was
mentioned in the introduction: it is the mechanism of alternation
between the  $q_n$ and $q'_n$ that is behind uniform
stretch for all $m \in \NN$ of the Birkhoff sums $S_m \varphi$ of an adequately chosen ceiling function $\varphi$. In addition to the first two, the third condition is useful to obtain rank one towers for the translation $R_{\a,\a'}$. It is easy to prove the existence of an uncountable and dense set of couples in ${\RR}^2$ satisfying 
 (\ref{10000})--(\ref{30000}) (see \cite{Y} and \cite{mls}).

\subsection{A special tower for $R_{\a,\a'}$.} \label{vvvspecialtower}

Let $(\a,\a') \in Y$ and consider on $\TT^2$ the translation $R_{\a,\a'}$ of vector $(\a,\a')$. 

\begin{defi} \label{R} For $0 \leq j \leq q_{n} q'_{n-1} -1$, define on ${\TT}^2$ the rectangles
$$R_{n}^j := \left( j { p_{n} \over q_{n}}, j { p_{n} \over q_{n}} + {1 \over q_n}  \right) \times \left( j {p'_{n-1}
\over  q'_{n-1}}, j {p'_{n-1}
\over  q'_{n-1}} + {1 \over q'_{n-1}}  \right).$$ 
\end{defi}

\begin{lemm} The
  rectangles $R_n^j$, $0\leq j \leq   q_{n} q'_{{n}-1} -1,$ are disjoint
 and their union is of full Haar measure on ${\TT}^2$.
\end{lemm}

\noindent {\sl Proof. }  Suppose $j$ and $j'$ are such that 
\begin{eqnarray*}
 \ j  \ {p_{n} \over q_{n}}  \  - \  j' {p_{n} \over q_{n}} \ \    & \in  \ZZ, \\
j {p'_{{n}-1} \over q'_{{n}-1}}   - j' {p'_{{n}-1} \over q'_{{n}-1}}  & \in  \ZZ. 
\end{eqnarray*}
\noindent Because $q_{n}$ and $p_{n}$ are relatively prime $q_{n}$ divides $j-j'$, and we have the
same for $q'_{{n}-1}$ but since we assumed that $q_{n}$ and $q'_{{n}-1}$ are relatively  prime $q_{n}
q'_{{n}-1}$ has to divide $j-j'$ and $j-j' \geq q_{n} q'_{{n}-1}$. Hence, up
to $ q_{n} q'_{{n}-1} -1 $ the  $R_n^j$ are indeed disjoint. \carre 

\vspace{0.2cm}

The rational translation  $R_{ \left( {p_{n} /  q_{n}}, \ {p'_{n-1}
/  q'_{n-1}} \right) }$ approximates the
    translation 
 $R_{\a,\a'}$ (equation (\ref{reduites1})), and the
 tower of the rational translation, $R_n^0$,..., $R_{n}^{{q_{n} q'_{n-1}-1}}$, is almost a tower for the
 irrational one. To this difference that the rectangle $R_{n}^0$
 is periodic under the action of $R_{ \left( {p_{n} /  q_{n}}, \ {p'_{{n}-1}
/  q'_{{n}-1}} \right) }$ while its first return on itself under the action of $R_{\a,\a'}$ is 
 shifted to the right on the $y$-axis by 
$\tn q_{n} q'_{{n}-1} \a'  \tn \sim {q_{n} / q'_{n}}$ (from (\ref{reduites2})). In the $x$ direction, the shift of the first return is far smaller since $\tn q_{n} q'_{{n}-1} \a \tn \sim {q'_{n-1} / q_{n+1}}$. This will allow us to select a special tower for $R_{\a,\a'}$ with base essentially the rectangle $[0, 1 /q_n] \times [0, q_n / q'_n]$.  
The stacking of the levels of the corresponding tower, from left to right in each
$R_n^j$, $0 \leq j \leq q_{n}q'_{{n}-1} -1$, displays a clear analogy with what happens
for an irrational rotation on the circle as well as in the cutting and stacking constructions and will be behind the cumulation of
staircase stretch by the Birkhoff sums of the ceiling function that we
will later consider over $R_{\a,\a'}$ (see Property (X.3) in Section \ref{pX}). 

\subsubsection{Description of the tower.}  We give now a precise description of the tower we want to consider.

\begin{defi} \label{B}   Denote by $[.]$ the integer part and let 
\begin{eqnarray*}  r_n := \left[ {q'_{n} \over q_{n} q'_{{n}-1}} \right] - 1, \end{eqnarray*}
and 
\begin{eqnarray*}  h_n := \left[ (1 - {2 \over n})r_n \right] q_{n} q'_{{n}-1}. \end{eqnarray*}
Define the rectangle
$${{B}}_n^0 := \left[{1 \over {n} q_{n}},
( 1 - { 1 \over {n}}){ 1 \over q_{n}} \right] \times  \left[{ 1 q_{n}\over {n}
q'_{n}}, (1-{1 \over {n}}){q_{n} \over q'_{n}} \right],$$
and denote by  ${{B}}_n^h$ its image under $h$ iterations of
$R_{{\a,\a'}}$. 
\end{defi}

One of the goals of this section is to prove the following
\begin{prop} \label{vvvtower} The translation  $R_{{\a,\a'}}$ is
rank one by rectangles: the sequence of towers ${\lbrace \GG (B_n^0, h_n) \rbrace}_{n \in \NN}$ is a rank one sequence for $R_{\a,\a'}$.
\end{prop}

\begin{rema} A more exhaustive tower over $B_n^0$ for $R_{\a,\a'}$ would be the one with  $r_n q_{n}q'_{n-1} \sim q'_n$ levels. The term $(1 - {2 \over n}) $ in the expression of $h_n$ is used to put aside the top levels over $B_n^0$ that will not lift to monochromatic levels for the flow but will rather carry uniform stretch (Compare Properties (Y.2) and (Y.3)-(Y.3') in Section \ref{pY}). 
\end{rema}

\begin{defi}(The rational rotation) \label{D} 
For every $0 \leq i \leq r_n$ define the disjoint subsets of $R_n^0$ 
\begin{eqnarray*}  D_{n}^{i q_{n} q'_{n-1} }  &:=& 
   \left[ {1 \over n^2} {1 \over q_{n}},   (1-{ 1 \over n^2}) {1 \over q_{n}} \right] \times \left[   i{q_{n} \over q'_{n}} + {q_{n} \over n^2
q'_{n}}, (i+1) {q_{n} \over q'_{n}} - { 1 \over n^2}
{q_{n} \over q'_{n}} \right] . \end{eqnarray*} 
For every $1 \leq j \leq q_n q'_{n-1} -1$ define the disjoint subsets of $R_n^j$
\begin{eqnarray*} D_{n}^{j+ i q_{n} q'_{n-1} }  &:=& R_{{p_n \over q_n}, {p'_{n-1} \over q'_{n-1}} }^j (D_{n}^{iq_nq'_{n-1}}). \end{eqnarray*}
\end{defi}
\begin{center}
\setlength{\unitlength}{0.00041667in}
\begingroup\makeatletter\ifx\SetFigFont\undefined%
\gdef\SetFigFont#1#2#3#4#5{%
  \reset@font\fontsize{#1}{#2pt}%
  \fontfamily{#3}\fontseries{#4}\fontshape{#5}%
  \selectfont}%
\fi\endgroup%
{\renewcommand{\dashlinestretch}{30}
\begin{picture}(16103,3999)(0,-10)
\path(1125,507)(11070,507)
\path(2205.000,4080.000)(2250.000,3960.000)(2295.000,4080.000)
\path(1575,3747)(1575,282)
\path(10755,3747)(10755,282)
\path(6165,3747)(6165,282)
\path(1080,1497)(11025,1497)
\path(1080,2442)(11025,2442)
\path(1665,1407)(2115,1407)(2115,597)
	(1665,597)(1665,1407)
\path(2250,1407)(2700,1407)(2700,597)
	(2250,597)(2250,1407)
\path(2835,1407)(3285,1407)(3285,597)
	(2835,597)(2835,1407)
\path(6255,3297)(6705,3297)(6705,2532)
	(6255,2532)(6255,3297)
\path(1080,3387)(11025,3387)
\path(6840,3297)(7290,3297)(7290,2532)
	(6840,2532)(6840,3297)
\path(7425,3297)(7875,3297)(7875,2532)
	(7425,2532)(7425,3297)
\path(5445,1407)(5895,1407)(5895,597)
	(5445,597)(5445,1407)
\path(10035,3297)(10485,3297)(10485,2532)
	(10035,2532)(10035,3297)
\put(5085,147){\makebox(0,0)[lb]{\smash{{{\SetFigFont{9}{10.8}{\familydefault}{\mddefault}{\updefault}     }}}}}
\put(10890,57){\makebox(0,0)[lb]{\smash{{{\SetFigFont{10}{12.0}{\familydefault}{\mddefault}{\updefault}$y$}}}}}
\put(1035,3792){\makebox(0,0)[lb]{\smash{{{\SetFigFont{10}{12.0}{\familydefault}{\mddefault}{\updefault}$x$}}}}}
\put(3690,867){\makebox(0,0)[lb]{\smash{{{\SetFigFont{9}{10.8}{\familydefault}{\mddefault}{\updefault} $R_n^0$}}}}}
\put(5310,57){\makebox(0,0)[lb]{\smash{{{\SetFigFont{9}{10.8}{\familydefault}{\mddefault}{\updefault}$D_n^{r_nq_nq'_{n-1}}$}}}}}
\put(1710,57){\makebox(0,0)[lb]{\smash{{{\SetFigFont{9}{10.8}{\familydefault}{\mddefault}{\updefault}$D_n^0$    \    $D_n^{q_nq'_{n-1}}...$}}}}}
\put(8280,2667){\makebox(0,0)[lb]{\smash{{{\SetFigFont{9}{10.8}{\familydefault}{\mddefault}{\updefault}$R_n^j$}}}}}
\put(6165,1902){\makebox(0,0)[lb]{\smash{{{\SetFigFont{9}{10.8}{\familydefault}{\mddefault}{\updefault}$D_n^j$   \  $D_n^{j+q_nq'_{n-1}}... \ \ \ \ \ \ \ \ \ \ \ \  \ $      $D_n^{j+r_nq_nq'_{n-1}}$}}}}}
\put(0,417){\makebox(0,0)[lb]{\smash{{{\SetFigFont{6}{7.2}{\familydefault}{\mddefault}{\updefault}.}}}}}
\end{picture}
}

\end{center}
$${\rm Fig.} \ 1. {\rm \ \ The \ rational \ rotation}. $$ 

\vspace{0.2cm}

\begin{defi} \label{RR} Let
 \begin{eqnarray*}  {\overline{R}}_n^0 := \left[{ 1 \over {n} q_{n}}, (1-{ 1
\over n}) {1 \over  q_{n}} \right] \times \left[{0}, (1- {1 \over n}) {1
\over  q'_{{n}-1}} \right] \subset R_{n}^0, \end{eqnarray*} 
and for every $ 0 \leq j \leq q_{n} q'_{{n}-1} -1$ let
$$\overline{R}_n^j:= R_{\a,\a'}^j (\overline{R}_n^0).$$
\end{defi}

The following Lemma encloses the facts that we will need about the combinatorics of $R_{\a,\a'}$ at step $n$. Recall that we have assumed that $n$ is odd so that  
 $ \a' - {p'_{{n}-1} \over q'_{{n}-1}} > 0 $.

\begin{lemm} \label{vvvzabb} For any $0 \leq j \leq  q_{n}q'_{{n}-1}-1 $, we have 
\begin{eqnarray} \label{hh1} {\overline{R}}_n^j \subset  R_n^j. \end{eqnarray}

\noindent Furthermore, \ for each  $0 \leq j \leq  q_{n} q'_{{n}-1} -1$, there is a real number $0 \leq \beta_{n,j} \leq  { q_{n} /
q_{n}'}$, such that for any $0 \leq i \leq r_n -1$, one has
\begin{eqnarray} \label{hh2} T_{0, - \beta_{n,j}} \left( B_n^{j +i q_{n} q'_{{n}-1}}
\right) \subset  {D}_{n,j+ i q_{n} q'_{n-1}}, \end{eqnarray}
where  $T_{0, \beta}$ is the translation of vector $(0, \beta)$.

\end{lemm} $ $

  An immediate corollary of the above lemma is that the sets $B_n^h$ 
are
disjoint for $ 0 \leq h \leq h_n$. Since in addition $h_n \lambda^{(2)} (B_n^0) \sim {( 1 -{2 \over n})}^3$ the Proposition \ref{vvvtower} follows. The above lemma also shows how to construct a function that is constant on the levels of the tower over $B_n^0$. 

\begin{coro} \label{vvvescalier} If a real function $\kappa$ defined on the torus is constant on $D_n^{j+ iq_{n} q'_{{n}-1}}$ for
some   $i \leq  r_n -1$ and some $ 0 \leq j \leq q_{n}
q'_{{n}-1} -1$, then the function 
$  \kappa \circ T_{0, - \beta_{n,j}}$ is constant on $ {B}_n^{j + i
  q_{n} q'_{{n}-1}}$.  In particular, if $\kappa$ is $(1 / q_n, 1 /q'_{n-1})$ periodic and constant on $D_{n}^{iq_nq'_{n-1}} \subset R_n^0$, then  for any $0 \leq j \leq q_n q'_{n-1} -1$ we have that $  \kappa \circ T_{0, - \beta_{n,j}}$ is constant on $ {B}_n^{j + i
  q_{n} q'_{{n}-1}}$. \end{coro}

$$$$
\begin{center}
\setlength{\unitlength}{0.00041667in}
\begingroup\makeatletter\ifx\SetFigFont\undefined%
\gdef\SetFigFont#1#2#3#4#5{%
  \reset@font\fontsize{#1}{#2pt}%
  \fontfamily{#3}\fontseries{#4}\fontshape{#5}%
  \selectfont}%
\fi\endgroup%
{\renewcommand{\dashlinestretch}{30}
\begin{picture}(11126,4494)(0,-10)
\put(0,4152){\makebox(0,0)[lb]{\smash{{{\SetFigFont{10}{12.0}{\familydefault}{\mddefault}{\updefault}$x$}}}}}
\put(10665,147){\makebox(0,0)[lb]{\smash{{{\SetFigFont{10}{12.0}{\familydefault}{\mddefault}{\updefault}$y$}}}}}
\path(180,597)(10710,597)
\path(540,4467)(540,57)
\dashline{60.000}(630,1542)(1125,1542)(1125,687)
	(630,687)(630,1542)
\path(180,1632)(10710,1632)
\dashline{60.000}(1260,1542)(1755,1542)(1755,687)
	(1260,687)(1260,1542)
\dashline{60.000}(1890,1542)(2385,1542)(2385,687)
	(1890,687)(1890,1542)
\path(135,2577)(10710,2577)
\path(540,1497)(4725,1497)(4725,732)
	(540,732)(540,1497)
\thicklines
\path(1305,1407)(1710,1407)(1710,822)
	(1305,822)(1305,1407)
\path(1935,1407)(2340,1407)(2340,822)
	(1935,822)(1935,1407)
\thinlines
\path(10215,4422)(10215,12)
\thicklines
\path(675,1407)(1080,1407)(1080,822)
	(675,822)(675,1407)
\path(6390,3387)(6795,3387)(6795,2802)
	(6390,2802)(6390,3387)
\thinlines
\path(5310,3927)(5535,3927)
\path(5265,4422)(5265,12)
\path(5310,3972)(5310,3882)
\path(180,3657)(10710,3657)
\dashline{60.000}(5985,3567)(6480,3567)(6480,2667)
	(5985,2667)(5985,3567)
\dashline{60.000}(6615,3567)(7110,3567)(7110,2667)
	(6615,2667)(6615,3567)
\path(5580,3477)(9765,3477)(9765,2712)
	(5580,2712)(5580,3477)
\dashline{60.000}(5355,3567)(5850,3567)(5850,2667)
	(5355,2667)(5355,3567)
\thicklines
\path(5760,3387)(6165,3387)(6165,2802)
	(5760,2802)(5760,3387)
\path(7020,3387)(7425,3387)(7425,2802)
	(7020,2802)(7020,3387)
\thinlines
\drawline(5085,2262)(5085,2262)
\path(5535,3972)(5535,3882)
\put(2880,1002){\makebox(0,0)[lb]{\smash{{{\SetFigFont{8}{9.6}{\familydefault}{\mddefault}{\updefault}$\overline{R}_n^0$}}}}}
\put(7785,2982){\makebox(0,0)[lb]{\smash{{{\SetFigFont{8}{9.6}{\familydefault}{\mddefault}{\updefault}  $\overline{R}_n^j$}}}}}
\put(675,147){\makebox(0,0)[lb]{\smash{{{\SetFigFont{8}{9.6}{\familydefault}{\mddefault}{\updefault}$B_n^0 \ \ B_n^{q_{2n}q'_{2n-1}}...$}}}}}
\put(5265,4107){\makebox(0,0)[lb]{\smash{{{\SetFigFont{8}{9.6}{\familydefault}{\mddefault}{\updefault} $\beta_{n,j}$}}}}}
\put(5760,2127){\makebox(0,0)[lb]{\smash{{{\SetFigFont{8}{9.6}{\familydefault}{\mddefault}{\updefault}$B_n^j \ \ B_n^{j+q_{2n}q'_{2n-1}}...$}}}}}
\end{picture}
}

\end{center}
$${\rm Fig.} \ 2. {\rm \ The \ special \ tower \ of \ } R_{\a,\a'}.  $$

\vspace{0.2cm}

\noindent {\sl Proof of Lemma \ref{vvvzabb}. } We will only prove (\ref{hh2}) that actually implies (\ref{hh1}) if we take $i=0$.  We will proceed separately for the $x$ and the $y$ direction.  Take a couple $(i,j)$ of integers such that $0 \leq j \leq  q_{n} q'_{{n}-1} -1 $ and $0 \leq i \leq  r_n$. Take $(x,y) \in  {B}_n^{j +i q_{n}
q'_{{n}-1}} $, hence
$$          {1 \over {n} q_{n}}  \leq   \lbrace {x} - j \a - iq_{n}
q'_{{n}-1} \a  \rbrace  \leq
    (1 -{ 1\over n}) {1 \over q_{n}}.  $$
From (\ref{reduites1}) we have 
$$ |  iq_{n} q'_{{n}-1} \a -iq'_{n-1}p_n | \leq i {q'_{n-1} \over q_{n+1}} \leq {q'_n \over q_n  q_{n+1}} = o({1 \over nq_n}),$$
and 
$$ \left| j \a - j {p_{n} \over q_{n}} \right| \leq {j \over q_{n+1}} \leq {q_n q'_{n-1} \over q_{n+1}} = o({1 \over nq_n}),$$
therefore 
\begin{eqnarray} \label{vvvxside} {1 \over n^2 q_{n}} \leq \lbrace x - j{p_{n} \over q_{n}} \rbrace \leq 
{1 \over q_{n}} - {1 \over n^2 q_{n}}. \end{eqnarray}

For the coordinate $y$ we have that  $(x,y) \in  {B}_n^{j +i q_{n}
q'_{{n}-1}} $ implies 
$$          {1 \over {n}} { q_{n} \over q'_{n}}  \leq   \lbrace {y} - j \a' - iq_{n}
q'_{{n}-1} \a'  \rbrace  \leq { q_{n} \over q'_{n}} - {1 \over {n}} { q_{n} \over q'_{n}}. $$
From (\ref{reduites2}) we have 
$$ \lbrace i q_{n} q'_{n-1} \a' - i { q_{n} \over q'_{n}}  \rbrace \leq 
 { i q_{n} q'_{n-1} \over q'_{n} q'_{n+1}} = o(  { q_{n} \over nq'_{n}}),$$
again from (\ref{reduites2}) we have 
$$ \left| j \a' - j {p'_{{n}-1} \over q'_{{n}-1}} - { j \over q'_{{n}-1} q'_{n} } \right| < { j
\over q'_{n} q'_{{n}+1}} = o( { q_{n} \over nq'_{n}}).$$ 
It follows that if we take 
 $ \beta_{n,j}$ to be   ${ j / q'_{{n}-1} q'_{n}} \leq {q_{n} / q'_{n}}$, we will have 
\begin{eqnarray} \label{vvvyside}
     {1 \over n^2}  { q_{n} \over q'_{n}}             \leq  \lbrace {y} - j {p'_{{n}-1} \over q'_{{n}-1}} - i { q_{n} \over q'_{n}} - \beta_{n,j}  \rbrace \leq  { q_{n} \over q'_{n}}  - {1 \over n^2}  { q_{n} \over q'_{n}}. \end{eqnarray}

The Lemma follows from (\ref{vvvxside}) and (\ref{vvvyside}). \carre

\section{The properties of  the function $ {X}_n$.}   \label{pX}

 The function $X_n$ will be  essentially a staircase function over every $R_n^j$ almost constant on  each level of the tower $\GG (B_n^0 , h_n )$ (see Section \ref{vvvspecialtower}). It will also be a trigonometric polynomial so that the derivatives of its Birkhoff sums over $R_{\a,\a'}$ must be uniformly bounded (in $m$). We will list here the properties required of $X_n$ and postpone its effective construction to Section \ref{chiconstruct}. Set

\begin{eqnarray} \label{epsn} \epsilon_n := {({q'_{n-1})}^7 \over  q'_{n}}. \end{eqnarray}

\begin{prop} \label{X} There exists a sequence of functions $X_n \in C^\infty(\TT^2, \RR)$ with the following properties,  

\begin{enumerate}

\item[${(X.1)}$] \label{chicontrolx} For any $r\geq 1$, there exists a constant $c(r)$ such that $\|D_x^r {X}_n \| \leq c(r) n^{2r +2} q_n^r \epsilon_n.$  
 
\vspace{0.1cm} 
 
 \item[${(X.2)}$] \label{chicontroly}  There exists a constant $c$ such that ${\|X_n\| + \|D_y {X}_n \|} \leq c n^4 {q'_n \over q_n} \epsilon_n.$
 
\vspace{0.1cm}

\item[${(X.3)}$] \label{chirankone} Let $B_n^0$ and $r_n$ be as in Definition  \ref{B}. For every $ j \leq q'_{n-1}q_{n} -1$ and
every $i \leq (1 -{2 \over n}) r_n$, we have for $(x,y) \in  B_n^{j + i q_{n}q'_{n-1}}$,  
$$|{X}_n (x,y) -  i \epsilon_n| \leq {1 \over  q_{n}q'_n}.$$

\vspace{0.1cm} 

\item[${(X.4)}$] \label{chitroncaturex} For $r \geq 1$, for $n$ large enough we have for every $m \in \NN$  
$$ \| { D_x^r  S_m \sum_{l \leq n}  {X}_l } \| \leq {q_{n+1} \over e^{q_n}}.  $$

\vspace{0.1cm}

  \item[${(X.5)}$] \label{chitroncaturey} For $p \geq 1$, for $n$ large enough we have for every $m \in \NN$ 
$$\| {D_y^p  S_m \sum_{l \leq n} {X}_l } \| \leq {q'_n}^{(p+4)}.$$  

\end{enumerate}

\end{prop} 

\begin{rema}[Choice of $\epsilon_n$] \label{racine} Due to (\ref{10000}) and (\ref{20000})  and our choice of $\epsilon_n$  in (\ref{epsn}) we have that (X.1) yields $\|D_x^r X_n \| \leq 1 / \sqrt{q'_n}$ and (X.2) yields $\|X_n\| + \|D_y X_n\| \leq 1 / \sqrt{q_n}$.   On the other hand, with this choice of $\epsilon_n$, (X.3) insures that already from $q_n / n^2$ the Birkhoff sums of $X_n$ are stretching in the $y$ direction above intervals of length $1 / {(q'_{n-1})}^6$ (see (\ref{sce})) while (X.5) implies that the lower order terms are almost constant  above such intervals.
\end{rema}

\section{The properties of the function $Y^{n}$.} \label{pY}

The function $Y_n$ is obtained from $\cos (2 \pi q_{n} x) / e^{q_{n}}$ multiplied by a    bump function essentially equal to 0
  over the rectangles  $\overline{R}_n^j \subset R_n^j$ of Definition \ref{RR} (i.e. over the first $h_n$ levels of the tower $\GG (B_n^0,h_n)$ and equal to 1 on the last $1 / n$ proportion of the rectangles $R_n^j$ (to produce uniform stretch when this end part of $R_n^j$ is visited). In addition $Y_n$ is taken to be a  trigonometric polynomial in order for  its Birkhoff sums over $R_{\a,\a'}$ to be uniformly bounded (in $m$). We will list here the properties required of $Y_n$ and postpone its effective construction to Section \ref{psiconstruct}.

\begin{prop} \label{Y} There exists a sequence of functions $Y_n \in C^\infty(\TT^2, \RR)$ with the following properties,  

\begin{enumerate}

\item[${(Y.1)}$] \label{psicontrol} For $r, p \in \NN$, there exists a constant $c'(r,p)$ such that we have 
$$\| D^r_x D^p_y {Y}_n  \| \leq c'(r,p) n^{3p} ({q'_{n-1})}^p {q_n^r \over e^{q_{n}}}. $$

\vspace{0.1cm}
 
\item[${(Y.2)}$]  \label{psirankone} For $h_n$ and $B_n^0$ as in Section \ref{vvvspecialtower}   we have   
\begin{eqnarray*}  \sup_{m \leq h_n} \sup_{z \in {B}_n^0} | S_m {{Y}}_n(z)| \leq{ 1  \over q'_{n}} . \end{eqnarray*}

\vspace{0.1cm}

\item[${(Y.3)}$] \label{psistretch}
For $n$ sufficiently large,  for  $m \in [q'_{n}, 2 q_{n+1} / {(n+1)}^2]$, for $x$ such that $\lbrace q_{n} x \rbrace  \in [{1 \over n}, {1 \over 2} -{1 \over n}] \cup  [{1 \over 2} + {1 \over n}, 1 -{1 \over n}]$ and for  any $y \in \TT$ we have
\begin{eqnarray*} 
\left| {  {D_x S_m {Y}}_n  }(x,y) \right| &\geq&
 {\pi \over 2 n^2}  { q_n \over
 e^{q_{n}}} m,     \\
\| {D^2_x  S_m {{Y}}_n} 
\| &\leq& 5 \pi^2 {q_n^2 \over
  e^{q_{n}}} m. \end{eqnarray*}

\vspace{0.1cm}

\item[${(Y.3')}$] 
For any $\eta>0$, we have for $n$ sufficiently large,  for  $m \in [q'_{n}/2n^2, q'_{n}]$,  for $x$ such that $\lbrace q_{n} x \rbrace  \in [{1 \over n}, {1 \over 2} -{1 \over n}] \cup  [{1 \over 2} + {1 \over n}, 1 -{1 \over n}]$ and for $y \in \TT$ such that $1-m/q'_n+\eta \leq \lbrace q'_{n-1} y \rbrace \leq 1 - \eta$  
\begin{eqnarray*} 
\left| {  {D_x S_m {Y}}_n  }(x,y) \right| &\geq&
 {\pi \over 2 n^2}  { q_n \over
 e^{q_{n}}} m,     \\
\| {D^2_x  S_m {{Y}}_n} 
\| &\leq& 5 \pi^2 {q_n^2 \over
  e^{q_{n}}} m. \end{eqnarray*}

\vspace{0.1cm}
    
\item[${(Y.4)}$] \label{psitroncaturex} For $r \geq 1$, for $n$ large enough we have for any $m \in \NN$ 
$$ \| { D_x^r S_m \sum_{l \leq n} {Y}_l } \| \leq  {q_{n+1} \over e^{{1 \over 2} q_n}}. $$
 
\vspace{0.1cm}

\item[${(Y.5)}$]  \label{psitroncaturey}  
 For $p \geq 1$, for $n$ large enough we have for any $m \in \NN$ 
$$\| {D_y^p S_m \sum_{l \leq n} {Y}_l  } \| \leq {q'_{n} \over e^{ {1 \over 2} q_n}}.$$ 
 
\end{enumerate} 

\end{prop}

\noindent {\sc Remark.} We stated (Y.3) and (Y.3') separately because they will
be used at different places in the proof, (Y.3) corresponding to a scale of time where uniform stretch in the $x$ direction is enough to yield mixing while (Y.3') is used at a scale of time where it yields mixing only on ``part'' of the space as will be explained in \S \ref{com}.

Define now 

$$ \varphi(x,y) := \varphi_0 + \sum_{n \geq n_0} X_n(x,y) + Y_n(x,y),$$
where $\varphi_0 \in \RR$ and $n_0$ are chosen such that $\varphi$ is strictly positive and has mean value one. From (X.1), (X.2) and (Y.1) we get that  $\varphi$ is of class $C^1$ on $\TT^2$ and is $C^\infty$ in the $x$ variable.  With $(\a,\a') \in Y$ and with the above properties on $X_n$ and $Y_n$ we will prove in the next two sections the following

\begin{theo} The special flow $T^t_{R_{\a,\a'},\varphi}$ is of rank one and mixing.
\end{theo}

\section{Proof of the rank one property.} \label{proof of rank one}

We want to check  Criterion \ref{criterion}  for $T^t_{R_{\a,\a'},\varphi}$. More precisely, given $h_n$ and $B_n^0$ as in Definition \ref{B} and having Proposition \ref{vvvtower},  we want to show that 

\begin{eqnarray} \label{rr}
\lim_{n \rightarrow \infty} \sup_{m\leq h_n} \sup_{z,z' \in B_n^0} |S_m \varphi(z') - S_m \varphi(z)| = 0.
  \end{eqnarray} 

Due to to the properties (\ref{10000}) and (\ref{20000}) of the sequences ${\lbrace q_n \rbrace}_{n \in \NN}$ and  ${\lbrace q'_n \rbrace}_{n \in \NN}$  we have for $n$ large enough:

\begin{itemize} 

\item (X.2) (Remark \ref{racine}) implies that  
$$\sup_{m\leq h_n} \|S_m (\sum_{l \geq n+1} X_l) \| \leq
 h_n \sum_{l \geq n+1} {1 \over \sqrt{q_l}} \leq  q'_n \sum_{l \geq n+1} {1 \over \sqrt{q_l}} \leq 
2 {q'_n \over \sqrt{q_{n+1}}},$$ 
 
\item (Y.1) implies that 
$$\sup_{m\leq h_n} \|S_m (\sum_{l \geq n+1} Y_l) \| \leq c'(0,0)q'_n \sum_{l \geq n+1} {1 \over e^{q_l}} \leq 2c'(0,0) {q'_n \over e^{q_{n+1}}},$$

\item (X.4) and (X.5) imply  that 
$$\sup_{m\leq h_n} \sup_{z,z' \in B_n^0} |S_m \sum_{l \leq n-1} X_l (z') - S_m \sum_{l \leq n-1} X_l (z)| \leq {1 \over e^{q_{n-1}}} +  {(q'_{n-1})}^5 { q_n \over q'_n},$$
 
\item (Y.4) and (Y.5) imply  that 
$$\sup_{m\leq h_n} \sup_{z,z' \in B_n^0} |S_m \sum_{l \leq n-1} Y_l (z') - S_m \sum_{l \leq n-1} Y_l (z)| \leq {1 \over e^{{1 \over 2} q_{n-1} }} + {q'_{n-1} \over e^{ {1 \over 2} q_{n-1}}} {q_n \over q'_n},$$
 
\item (X.3) implies that 
$$\sup_{m\leq h_n} \sup_{z,z' \in B_n^0} |S_m X_n (z') - S_m X_n (z)| \leq {2 \over q_n}.$$

\end{itemize}

Together with (Y.2) the above estimations yield the required (\ref{rr}).

$ \ \ $  \carre

\section{Proof of mixing.}  \label{proof of mixing}

We will prove mixing in three steps depending on the range of $t \in \RR$. In step 1, mixing is obtained for some range of time due to uniform stretch of the Birkhoff sums of $Y_n$ (Property (Y.3)). In step 2 mixing is obtained for another range of time due to staircase stretch of the Birkhoff sums of $X_n$ (Property (X.3)). For the remaining times mixing is established in step 3 due to a combination of uniform stretch and staircase stretch mechanisms. The proof of mixing in steps 2 and 3 uses the existence of mixing intervals of time established in step 1.

\subsection {Step 1. Uniform stretch.} We will prove in this step that the sequence of intervals 
$[2q'_n, {q_{n+1} / {(n+1)}^2}]$ is mixing for the special flow $T^t_{R_{\a,\a'},\varphi}$ as in Definition \ref{m132}. 

\begin{defi} \label{P} Let ${\Omega}_n$ be a partition of the set 
$$J_n:= \left\{ x \in \TT \ / \  \lbrace q_{n} x \rbrace  \in [{1 \over n}, {1 \over 2} -{1 \over n}] \cup  [{1 \over 2} + {1 \over n}, 1 -{1 \over n}]  \right\}$$
in intervals of length between ${1 \over 2} e^{-q_n}$ and $e^{-q_n}$. Clearly ${\Omega}_n$ converges to the partition into points of $\TT$ as $n \rightarrow 0$. 
\end{defi}

We want to apply Criterion \ref{criterion mixing} to  the sequence ${\lbrace \Omega_n \rbrace }_{n \in \NN}$ and show that any sequence $t_n \in [2q'_n, {q_{n+1} / {(n+1)}^2}]$ is mixing.
Since $\varphi$ is continuous and has mean value $1$ we have by unique ergodicity of $R_{\a,\a'}$ that for sufficiently large $n$, for any $t \in 
[2q'_n, {q_{n+1} / {(n+1)}^2}]$, for any $(x,y) \in \TT^2$ 
$$m(x,y,t) \in [q'_n, 2{q_{n+1} / {(n+1)}^2}]$$
where $m(x,y,t)$ is as in Definition \ref{mzt}. Hence Step 1 will follow from Criterion \ref{criterion mixing} if we prove that for any interval $I \in {\Omega}_n$ and any $y \in \TT$ we have for  $m \in [q'_{n}, 2 q_{n+1} / {(n+1)}^2]$
\begin{eqnarray} \label{M} S_m \varphi (.,y) {\rm \ is \ } (K_n, \nu_n)-{ \rm \ uniformly \ stretching \ on \ }  I \times \lbrace y \rbrace, \end{eqnarray} 
for some sequences $K_n \rightarrow \infty$ and $\nu_n \rightarrow 0$.

To get this we deduce from the properties of $X_n$ and $Y_n$ the following estimates for $i=1,2$ and for $n$ large enough:

\begin{itemize} 
 
\item  (X.4) and (Y.4) imply that 
$$\mathop{\sup} \limits_{m \in \NN}  \| D_x^i S_m  \sum_{ l \leq n-1 } ( X_l+Y_l) \| \leq  2 {q_n e^{ -{1 \over 2} q_{n-1}}}, $$

\item (X.1) implies that $  \| D_x^i S_m (\sum_{l \geq n} X_l ) \| \leq m {2 \over \sqrt{q'_n}}$,

\item (Y.1) implies that $ \displaystyle \|D_x^i S_m (\sum_{l \geq n+1} Y_l ) \| \leq m {1 \over e^{{1 \over 2} q_{n+1}}}$.

\end{itemize}

The above quantities being negligible with respect to $m q_n / n^2 e^{q_n}$ if $m \in [q'_n, 2 q_{n+1} / {(n+1)}^2],$  we deduce from (Y.3) that for this range of $m$ and for $n$ large enough we have for $x \in J_n$, $y \in \TT$ 
 
\begin{eqnarray*} 
\left| {  {D_x S_m {\varphi}}  }(x,y) \right| &\geq&
 {\pi \over 3 n^2}  { q_n \over
 e^{q_{n}}} m,     \\
\| {D^2_x  S_m {{\varphi}}} 
\| &\leq& 6 \pi^2 {q_n^2 \over
  e^{q_{n}}} m. \end{eqnarray*}

Hence (\ref{10000}) implies that for any interval $I \in {\Omega}_n$ (Definition \ref{P}) we have 
for $m$ as above and any $y \in \TT$
\begin{eqnarray*} 
\inf_{x \in I} |D_x S_m \varphi (x,y)| |I| &\geq& {\pi q_n \over 6 n^2} {m \over e^{2q_n}} \\
&\geq& {\pi q_n \over 6 n^2} e^{q_n}, 
\end{eqnarray*} 
and since $|I| \leq  e^{-q_n}$ we get 
\begin{eqnarray*} 
\| D_x^2 S_m \varphi \| |I| \leq  {18 \pi n^2 q_n \over e^{q_n}}  \inf_{x \in I} |D_x S_m \varphi (x,y)|, \end{eqnarray*}
and the desired (\ref{M}) follows from Lemma \ref{deriveeseconde} with $K_n = {\pi q_n \over 6 n^2} e^{q_n}$ and $\nu_n =   {18 \pi n^2 q_n \over e^{q_n}} $.   \carre

\subsection{Step 2. Staircase stretch.} 
  We will prove in this step that the sequence of intervals 
$[q_n / n^2, {q'_n / n^2}]$  is mixing for the special flow $T^t_{R_{\a,\a'},\varphi}$ as in Definition \ref{m132}. From now on we will assume that $A$ is a fixed ball in $M_{R_{\a,\a'},\varphi}$.

 \subsubsection{ Consequence of Step 1.} We begin with a preliminary Lemma that is due to the existence of mixing intervals obtained in Step 1 and that will be useful in establishing mixing in this Step 2. 
\begin{lemm} \label{51} There exists a sequence of positive numbers  $\varepsilon_{1,l} \rightarrow 0$ such that if $\delta >0$ and $H \in \NN$ satisfy $\delta H \in [3 (l+1) {q'_l}, q_{l+1} / {(l+1)}^2]$ and $H \geq l^2$ then there exists a set ${\mathcal U} (H, \delta) \subset M_{R_{\a,\a'}, \varphi}$ with $\mu ({\mathcal U}) \geq 1- \varepsilon_{1,l}$  such that  for any $z \in {\mathcal U}$ we have 
\begin{eqnarray*}  \left|  {1 \over H} \sum_{i=0}^{H-1} \chi_A (T^{- i \delta }z) - \mu (A)  \right| < \varepsilon_{1,l}
\end{eqnarray*}  
where $C$ is a constant. 
\end{lemm} 

\noindent {\sl Proof. }  Given $f, g \in L^2 (M_{R_{\a,\a'}, \varphi}, \RR)$, we define the scalar product 
$$<f / g> := \int_{M_{R_{\a,\a'}, \varphi}} f(z) g(z) d \mu (z).$$
The lemma clearly follows if we prove that there exists a sequence $\varepsilon_l \rightarrow 0$ such that for any $f \in L^2 ( M_{R_{\a,\a'}, \varphi}, \RR)$ with $\int_{M_{R_{\a,\a'}, \varphi}} f(z) d\mu (z) =0$ we have 
 \begin{eqnarray} \label{pppp} 
 { \left\| {1 \over H} \sum_{i=0}^{H-1 } f \circ T^{- i \delta}  \right\|}_{L^2}  \leq C(f) \varepsilon_l.  \end{eqnarray}

It follows from Step 1 that there exists $\theta_l \rightarrow 0$ such that for any $\tau \in [2 q'_{l}, q_{l+1} / {(l+1)}^2]$ we have   
$$|<f \circ T^{-\tau} / f > | \leq C(f) \theta_l.$$
Hence we define  
\begin{eqnarray*} 
H' &=& \left[ {H \over l+1} \right], \\
\tau &=&  H' \delta,
\end{eqnarray*}
and we see that since $\tau \in [2 {q'_l}, q_{l+1} / {(l+1)}^3]$, we have 
for any $1 \leq j \leq l $, $|<f \circ T^{-j\tau} / f > | \leq C(f) \theta_l,$ then
$$ {\left( {\left\| \sum_{j=0}^{l}  f \circ T^{-j\tau}  \right\|}_{L^2} \right)}^2 \leq (l+1) {\| f \|}_{L^2}^2 + (l^2+l) C(f) \theta_l,$$
which gives 
 $$ {\left\| {1 \over l+1} \sum_{j=0}^{l}  f \circ T^{-j\tau}  \right\|}_{L^2} \leq {1 \over \sqrt{l+1}} {\| f \|}_{L^2} +  {C(f)}^{1 \over 2} \theta_l^{1 \over 2},$$ 
and since the measure $\mu$ is invariant by the flow we have
\begin{eqnarray*} 
 { \left\| {1 \over H'(l+1)} \sum_{i=0}^{H'(l+1)-1 } f \circ T^{- i \delta}  \right\|}_{L^2} &=& {\left\| {1 \over H'(l+1)} \sum_{p=0}^{H'-1} \sum_{j=0}^{l}  f \circ T^{-j\tau - p \delta}  \right\|}_{L^2} \\
 &\leq&  {1 \over \sqrt{l+1}} {\| f \|}_{L^2} +  {C(f)}^{1 \over 2} \theta_l^{1 \over 2}. \end{eqnarray*} 
We conclude  observing that \begin{eqnarray*}  { \left\| {1 \over H} \sum_{i=0}^{H-1 } f \circ T^{- i \delta} \right\| }_{L^2} \leq  { \left\| {1 \over H'(l+1)} \sum_{i=0}^{H'(l+1)-1 }   f  \circ T^{- i \delta} \right\| }_{L^2}  + {l+1 \over H} {\| f \|}_{L^2} \end{eqnarray*}
and  using the hypothesis $H \geq l^2$. \carre

\subsubsection{Good partial partitions at time $t$.} \label{col}  Recall from Section \ref{good} that to prove that 
$[q_n / n^2, q'_n / n^2]$            is an interval of mixing it is enough to show that for any $t \in [q_n / n^2, q'_n / n^2]$ there exists a sequence of partial partitions  $\Omega_t$ with sets of the form $\xi = R \times \lbrace s \rbrace$, $R \subset \TT^2$, such that  $\Omega_t$ converges to the partition into points of 
$M_{R_{\a,\a'}, \varphi}$ as $t \rightarrow \infty$ and such that for any $\varepsilon > 0$ we have for $n$ large enough  (\ref{fub2})  for any set ${\xi} \in \Omega_t$, that is
                                                
\begin{eqnarray*} 
\left| \lambda^{(2)} (\xi  \cap T^{-t_n} A) - \lambda^{(2)} (\xi)  \mu ( A ) \right| \leq  \varepsilon \lambda^{(2)}  (\xi),  
\end{eqnarray*}

In all this section we will assume $t \in [q_n / n^2, q'_n / n^2]$ is fixed. In relation with Lemma  \ref{51} we give the following 
 
\begin{defi} \label{5H} Given $\epsilon_n$ as in (\ref{epsn}) let 
\begin{eqnarray*} 
M &:=& [t] \\
\delta &:=& M \epsilon_n \\
 H &:=& \begin{cases}

\left[ {{q'_{n-2}}^8 \over \delta} \right]  \text{if $ M \leq e^{q_n}$} \\ 
\\
  \left[ {{(q'_{n-1})}^8 \over \delta} \right]  \text{if  $M > e^{q_n}$.} \end{cases}
\end{eqnarray*}
\end{defi} 

\begin{rema} \label{5remark} It is easy to see that in both cases we have $ H \geq n^2 $ and hence that we can apply Lemma \ref{51} to the couple $(\delta,H)$ with $l=n-2$ if $M \leq e^{q_n}$ and $l=n-1$ if $M > e^{q_n}$. It is also clear that $H \leq {q'_n \over q_n {(q'_{n-1})}^6}$. The latter will be crucial when we will want to prove that, for $m$ comparable to $t$ and above the sets of  "length" $H$ (Definition \ref{5subdivision}),  only the Birkhoff sums of $S_m X_n$ are responsible for the variations of $S_m \varphi$ (see Lemma \ref{57}).

\end{rema} 

\begin{defi} \label{5subdivision} With the notations of Definitions \ref{B} and \ref{5H}, we will call a set ${\xi} \times \lbrace s \rbrace  \subset M_{R_{\a,\a'}, \varphi}$ (${\xi} \subset \TT^2$ and $0 \leq s \leq  \inf_{(x,y) \in \xi} \varphi(x,y)$) {\sl good} if  
\begin{eqnarray} \label{5R} {\xi} :=  B_n^{j + i_0 q_{n} q'_{n-1}}  \cup  B_n^{j + (i_0 +1) q_{n} q'_{n-1}} \cup... B_n^{j + (i_0 +H -1) q_{n} q'_{n-1} } \end{eqnarray} 
where $0 \leq j \leq q_n q'_{n-1} -1$ and $i_0 \in \NN$ satisfies 
\begin{eqnarray} \label{5C1} i_0 + H < (1 -{ 4 \over n}) r_n \end{eqnarray}
and if there exists a point $z_0 \in B_n^{j + i_0 q_n q'_{n-1}} \times \lbrace s \rbrace $ such that  
\begin{eqnarray} \label{5C2}   z_0 \in T^{-t} {\mathcal U} \end{eqnarray}
where ${\mathcal U} = {\mathcal U}(\delta, H)$ is the set obtained in Lemma \ref{51} (with $l=n-2$ if  $M \leq e^{q_n}$ and $l=n-1$ if $M > e^{q_n}$). 
  
\end{defi}

From the fact that for $n$ large $H$ is negligible with respect to $r_n$ (see Remark \ref{5remark}) and the fact that  the measure of ${\mathcal U}$ can be made arbitrarily close to 0, and from what was said in \S \ref{col} we will finish if we prove (\ref{fub2}) for   any good set ${\xi} \times \lbrace s \rbrace$.

A set ${\xi} \times \lbrace s \rbrace $ being given we denote it for simplicity by ${\xi}$ and denote the sets $B_n^{j + (i_0 + i) q_{n} q'_{n-1}}$ by $B(i)$ for $0 \leq i \leq H -1$. We can also assume that $s = 0$ since this does not alter the proof.

Fix hereafter $\varepsilon > 0$. Fix also two balls in $M_{R_{\a,\a'}, \varphi}$, $A^+$ and $A^-$  such that $A \subset {\rm int} (A^+)$  and  $A^- \subset  {\rm int} (A)$ and such that 
\begin{eqnarray} \label{5eps} (1- \varepsilon^2) \mu(A^+) \leq \mu(A) \leq (1 + \varepsilon^2) \mu(A^-). \end{eqnarray}

The following Proposition encloses the essential consequence of  the staircase stretch displayed by the Birkhoff sums of $\varphi$ due to our definition of $X_n$ and $Y_n$.

\begin{prop} \label{5A} There exists $n_0 \in \NN$ such that given any    time $t  \in [q_n / n^2, q'_n / n^2]$ and any good   set ${\xi}$ (see Definition \ref{5subdivision}), we have for any $0 \leq i \leq H-1$ and any $z_0 \in B(0)$

\begin{itemize} 

\item If $T^{t - i \delta} (z_0) \in A^-$ then $ T^{t} (B(i)) \subset A$.

\item If $ T^{t} (B(i)) \cap A \neq  \emptyset $ then  $T^{t - i \delta} (z_0) \in A^+$.
\end{itemize}
\end{prop}

Before proving this Proposition we show how to derive (\ref{fub2}) from it. For $n \geq n_0$ and $z_0 \in B(0)$ we have

\begin{eqnarray*} \lefteqn{\sum_{i= 0}^{H-1} \chi_{A^-} (T^{t - i \delta}  (z_0)) \lambda^{(2)} (B(0)) \leq \lambda^{(2)} \left( {\xi} \cap T^{-t} A     \right) } \hspace{2cm}  \\ 
& & \leq  \sum_{i=0}^{H-1} \chi_{A^+} ( T^{t - i \delta} (z_0)) \lambda^{(2)}  (B(0)). \end{eqnarray*}
Lemma \ref{51} is applicable due to Remark \ref{5remark}, hence considering the latter equation for $z_0 \in B(0) \cap T^{-t} { \mathcal U}$ (see (\ref{5C2})) we obtain if $\varepsilon_{1,n-2}$ and $\varepsilon_{1,n-1}$ are sufficiently small that 
$$ 
 (1 - \varepsilon^2) H \mu(A^-)  \lambda^{(2)} (B(0)) 
 \leq \lambda^{(2)} \left( {\xi}  \cap T^{-t} A     \right) \leq  (1+  \varepsilon^2) H \mu ( A^+ )  \lambda^{(2)} (B(0)). $$
Since $H \lambda^{(2)} (B(0)) = \lambda^{(2)} ({\xi})$ this last inequality and (\ref{5eps}) lead to (\ref{fub2}) if $\varepsilon \leq {1 \over 4}$. 

\carre

\vspace{0.2cm}

In our proof of Proposition \ref{5A} we will need the following Lemma

\begin{lemm} \label{57} There exists a sequence $\varepsilon_{2,n} \rightarrow 0$ such that  if $0 \leq i_1 \leq i_2 \leq (1 -{ 4 \over n}) r_n $ and $i_2 - i_1 \leq q'_n / (q_n {(q'_{n-1})}^6) $  and $m \leq 2q'_n / n^2$ then for any $0 \leq j \leq q_n q'_{n-1}-1$ we have for any $z_{1} \in B^{j + i_1 q_n q'_{n-1} }$ and  $z_{2} \in B^{j + i_2 q_n q'_{n-1}}$ 
\begin{eqnarray} \label{sce}
\left| S_m \varphi(z_{2}) - S_m \varphi(z_{1}) - (i_2-i_1) m \epsilon_n \right| \leq \varepsilon_{2,n} \end{eqnarray}
 \end{lemm} 

\noindent {\sl Proof.} For $j, i_1, i_2,$ and $ m$ as above we have for every $l \leq m$ that $j + i_1 q_n q'_{n-1} +l  \leq  j + i_2 q_n q'_{n-1} +l \leq (1 - {3 \over n}) r_n q_n q'_{n-1}$ hence (X.3) implies that 
\begin{eqnarray*}
\left| X_n (R^l_{\a,\a'} z_{2}) -  X_n ( R^l_{\a,\a'} z_{1}) - (i_2-i_1) \epsilon_n \right| \leq {1 \over q_n q'_n} \end{eqnarray*}
hence 
\begin{eqnarray}
\left| S_mX_n (z_{2}) - S_m X_n ( z_{1}) - (i_2-i_1) m \epsilon_n \right| \leq {2 \over n^2 q_n}.  \end{eqnarray}

We still have to bound $|S_m (\varphi - X_n)(z_2) - S_m (\varphi - X_n)(z_1)|$. The condition $i_2 - i_1 \leq  q'_n / (q_n {(q'_{n-1})}^6) $ implies that the distance between the $y$ coordinates of $z_2$ and $z_1$ is less than $1 / {(q'_{n-1})}^6$ (see Definition \ref{B}). Therefore (X.4) and (X.5) imply that  $|S_m \sum_{l \leq n-1} X_l (z_2) -  
S_m \sum_{l \leq n-1} X_l (z_1)| \leq 1 / e^{q_{n-1}} + 1 / q'_{n-1}$. 

The control of the other terms in $S_m \varphi$ is similar to the one obtained in the  proof of the rank one property in Section \ref{proof of rank one}. \carre

\subsubsection{ Proof of Proposition \ref{5A}.}  We will prove the first point in the proposition, the second one being obtained similarly.  Let $0 \leq i \leq H-1$  and denote by $z_i$ some arbitrarily fixed point in $B (i)$. Define $ V \in \RR$ by 

\begin{eqnarray}
V  := t - S_{M} \varphi(z_i).  \label{5V}
  \end{eqnarray}

We will need in the sequel an upper bound on $|V|$:
\begin{lemm} \label{Vbound}  For $n$ sufficiently large, we have for any $m \in \NN$ and any $z \in \TT^2$ 
$$\left| S_m \varphi(z) - m \right| \leq 6 {m \over \sqrt{q_n}} + 2 {q_n \over e^{{1 \over 2} q_{n-1}}} + 2 {(q'_{n-1})}^5.$$
In particular since $M = [t]$ we have 
 \begin{itemize} 
\item If ${q_n \over n^2} \leq  M \leq e^{q_n}$ then $|V| \leq {M \over e^{q'_{n-2}}},$
\item If $e^{q_n} \leq M \leq {q'_n \over n^2}$ then $|V| \leq {M \over e^{q'_{n-1}}}$.
\end{itemize}
\end{lemm} 

We will prove this Lemma at the end of the section.

Define now
\begin{eqnarray} 
\label{5U} U   := t - i \delta - S_{M} \varphi(z_0). \end{eqnarray}
By definition of a special flow we have from (\ref{5U}) and (\ref{5V}) 
\begin{eqnarray*} 
T^{t - i \delta } (z_0) &=& T^{U  }  \left( R_{\a,\a'}^{M} z_0 \right), \\
T^{t  } (z_i) &=& T^{V}  \left ( R_{\a,\a'}^{M} z_i \right).
\end{eqnarray*}

But from Lemma \ref{57} and Remark \ref{5remark} it follows that as $n$ goes to $\infty$  
\begin{eqnarray*} \label{wa7ad} |U - V| \rightarrow 0, \end{eqnarray*}
hence if we consider a ball ${A^-}'$ strictly included in $A$ and strictly including $A^-$ (i.e. $ {A^-}'  \subset {\rm int} A$, $ A^- \subset {\rm int} {A^-}'$) we have for sufficiently large $n$ that  if as in the statement of the Proposition $T^U (R^M_{\a,\a'} z_0) = T^{t - i \delta} (z_0) \in A^-$ then 
\begin{eqnarray} \label{5vv}  T^{V} \left( R_{\a,\a'}^{M} z_0 \right) \in {A^-}'. \end{eqnarray}
To finish we must prove that for sufficiently large $n$ the latter implies that  
\begin{eqnarray}  
\label{5vvv}  T^{V} \left( R_{\a,\a'}^{M} z_i \right) \in {A}.
\end{eqnarray}
For this it is enough to show that as $n$ goes to infinity
\begin{eqnarray} \label{kj}   \sup_{m \leq 2|V|} \left| S_m \varphi (R_{\a,\a'}^M z_i) - S_m \varphi (R_{\a,\a'}^M z_0 )  \right|  \rightarrow 0.
\end{eqnarray}
Since  $2|V| + M \leq 2q'_n / n^2$ we have by Lemma \ref{57} again
 that for any $m \leq 2 |V|$ 
$$\left| S_m \varphi (R_{\a,\a'}^M z_i) - S_m \varphi (R_{\a,\a'}^M z_0 ) - im\epsilon_n \right| \leq \varepsilon_{3,n},$$ 
hence to get (\ref{kj}) and finish we just have to check that 
\begin{eqnarray*} |V| H \epsilon_n \rightarrow 0.
\end{eqnarray*} 
In light of the Definition \ref{5H} of $H$ and Lemma \ref{Vbound} we have two cases: If   $M \leq e^{q_n}$, then $|V| H \epsilon_n = |V| {q'_{n-2}}^8 / M \leq   {q'_{n-2}}^8 / e^{q'_{n-2}}$; If $M > e^{q_n}$ then  $|V| H \epsilon_n =  
 |V| {(q'_{n-1})}^8 / M \leq   {(q'_{n-1})}^8  / e^{q'_{n-1}}$.

It only remains to give the 

\noindent{ \sl Proof of  Lemma \ref{Vbound}.}  As in the proof of the rank one property in Section \ref{proof of rank one}, it follows from Propositions \ref{X} and \ref{Y} that  for $n$ sufficiently large we have for any $m \in \NN$ 
$$\left\| S_m  \sum_{l \geq n} (X_l + Y_l)  \right\| \leq 3{m \over \sqrt{q_n}},$$
while   
$$\left\| D_x S_m  \sum_{l \leq n-1}  (X_l + Y_l) \right\| \leq 2 {q_n \over e^{{1 \over 2} q_{n-1}}},$$
and  
$$ \left\| D_y S_m  \sum_{l \leq n-1}  (X_l + Y_l) \right\| \leq 2{(q'_{n-1})}^5,$$
which yields for any $z, z' \in \TT$
$$\left| S_m \varphi (z) - S_m \varphi (z') \right|  \leq 6 {m \over \sqrt{q_n} } +  2 {q_n \over e^{{1 \over 2} q_{n-1}}} + 2 {(q'_{n-1})}^5,$$
integrating along $z'$ we get 
$$\left| S_m \varphi (z) - m  \right|  \leq 6 {m \over \sqrt{q_n} } +  2 {q_n \over e^{{1 \over 2} q_{n-1}}} + 2 {(q'_{n-1})}^5,$$
from which Lemma \ref{Vbound} easily follows due to the inequalities (\ref{10000}) and (\ref{20000}) between the denominators of the best approximations of $\a$ and $\a'$.  \carre

\subsection{Step 3. Combining  staircase stretch and uniform stretch.} \label{com} In this step we want to complete the proof of mixing by showing that the intervals $[q'_n / n^2, 2q'_n]$  form a sequence of mixing intervals of time for the special flow. In this range of time, both mechanisms of mixing displayed in Step 1 and Step 2 enter into play and imply mixing for sets lying in different parts of $M_{R_{\a,\a'}, \varphi}$. In all this section we assume $t \in [q'_n / n^2, 2q'_n]$ is fixed and introduce 
\begin{eqnarray}  \theta := {t \over  q'_n}. \label{the} \end{eqnarray}
With the definition of $m(z,t)$ given in Section \ref{special flow}, observe that by unique ergodicity of $R_{\a,\a'}$ and continuity of $\varphi$ and since we chose $\int_{\TT^2} \varphi(x,y) dx dy =1$ then for $\eta > 0$ arbitrarily small there exists $t_0 \geq 0$ such that for any $t \geq t_0$ and any $z \in M_{R_{\a,\a'}, \varphi}$ we have 
\begin{eqnarray} \label{mz} (1 -\eta^2) t \leq m(z,t) \leq (1 + \eta^2) t. \end{eqnarray}  
Introduce the subsets of $M_{R_{\a,\a'}, \varphi}$ corresponding to uniform stretch and staircase stretch respectively 
\begin{eqnarray*} 
M^u (\theta, \eta)&:=& \left\{ (x,y,s) \in M_{R_{\a,\a'}, \varphi} \ / \   1 - \theta + \eta \leq \lbrace q'_{n-1}y \rbrace \leq 1 - \eta \right\} \\
M^s(\theta, \eta) &:=& \left\{ (x,y,s) \in M_{R_{\a,\a'}, \varphi} \ / \ \lbrace q'_{n-1}y \rbrace \leq 1 - \theta - \eta  \right\},  
\end{eqnarray*}
where for instance if $\theta \geq 1- \eta$ we have $M^s = \emptyset$. 

We see now how the combination of uniform and staircase stretch occurs:

\noindent $\bullet$ Due to (Y.3') we can repeat exactly the arguments of Step 1 and define in $M^u(\theta, \eta)$ a collection ${\Omega}^u (t)$ consisting  of intervals in the $x$ direction as in Definition \ref{P}
covering all but an arbitrarily small proportion  (as $n \rightarrow \infty$) of $M^u(\theta, \eta)$  
for which due to uniform stretch (\ref{fub2}) holds at time $t$ (for an arbitrary ball $A$ and an arbitrary precision $\varepsilon$ provided $n$ is large enough).

\vspace{0.2cm}

\noindent $\bullet$  For any $z \in M^s(\theta, \eta)$,  for any $u \leq t$ we have  $T^u(z) \in M^s(0, {\eta \over 2})$ (follows from the arithmetics of $\a$) hence all the calculations of Step 2 are still valid at this time $t$  for  the points in $M^s(\theta, \eta)$. Hence we can consider a collection $\Omega^s (t)$ consisting   of good sets as in Definition \ref{5subdivision}, covering all but an arbitrarily small measure of  $M^s(\theta, \eta)$, for which (\ref{fub2}) holds  (for an arbitrary ball $A$ and an arbitrary precision $\varepsilon$ provided $n$ is large enough).

Since the measure of the set of $M_{R_{\a,\a'}, \varphi}$ that is not included in 
\newline  $M^s(\theta, \eta) \cup M^u(\theta, \eta)$
 converges to $0$ as $\eta \rightarrow 0$ the partial partition ${\Omega}^u(t)  \cup \Omega^s (t)$ converges to the partition into points of $M_{R_{\a,\a'}, \varphi}$ as $\eta \rightarrow 0$ and $n \rightarrow \infty$ and Step 3 is thus completed.    \carre

\section{\sc Construction of $X_n$.} \label{chiconstruct}

\subsection{Construction of a first function $\tilde{X}_n$.}

Let $\theta$ be a $C^{\infty}$ increasing function on ${\RR}$ such that,  
$$ \theta(s) = 0 { \rm \ \ \ \ for \ \  \ \ } s \leq 0,$$
$$ \theta(s) = 1 { \rm \ \ \ \ for \ \  \ \ } s \geq 1.$$

With $r_n$ as in Section \ref{vvvspecialtower} and  $\epsilon_n$ as in (\ref{epsn}), 
consider on
$[0, { 1 \over q'_{{n}-1}}]$ the following function

\begin{eqnarray} \label{bas} \kappa_n  (y) = \left( \sum_{l= 1 }^{r_n}
\epsilon_n \theta \left(n^2 {q'_{n} \over q_{n}} (y - l {q_{n} \over q'_{n}})
\right) \right) \left(1 - \theta(nq'_{n-1}y-n+2) \right)  
\end{eqnarray}
extended to a $C^{\infty}$ function on ${\TT}^2$ independent of the
variable $x$ and of period ${ 1 / q'_{{n}-1}}$ in the variable
$y$. This is  possible because the left hand side in the above
expression is identically zero when $y \leq {q_{n} / q'_{n}},   $ 
while the right hand  side is identically zero for $ y \geq (1 - { 1 \over n}) {1
 / q'_{{n}-1}}.$

In relation with (X.1), (X.2) and (X.3) we will need the following Lemmas on $\kappa_n$

\begin{lemm} \label{kkkkk} With $r_n$ and $D_{n}^{i q_{n} q'_{{n}-1}} $ as in Definitions \ref{B} and \ref{D} we have for every $i \leq   (1 - { 2 \over n} ) r_n  $ that
${\kappa}_{n}$ is constant on  $D_{n}^{iq_{n}q'_{{n}-1}}$.  More precisely 
$$ {\kappa_n}_{ | D_{n}^{i q_{n}q'_{n-1}}} \equiv  i \epsilon_n. $$
\end{lemm}

\noindent {\sl Proof. } When $i \leq   (1 - { 2 \over n} ) r_n  $,  $(x,y) \in  D_{n}^{iq_{n}
q'_{{n}-1}}$ implies that $0 \leq y \leq (1 - { 2 \over n}) {1 \over
q'_{{n}-1}}$, in this case the right hand side in the expression of $\kappa_n$ is equal to 1 and 
$$ \kappa_n(y) = \sum_{l= 1}^{r_n}
\epsilon_n \theta \left(n^2 {q'_{n} \over q_{n}} (y - l {q_{n} \over q'_{n}})
\right),$$
but for every $l$, $\theta \left( n^2 {q'_{n} /  q_{n}} (y - l {q_{n}
/  q'_{n}}) \right)$  is constant equal to $0$ when $y \leq l {q_{n}
 / q'_{n}}$, and constant equal to $1$ when $ y \geq l {q_{n} / q'_{n}} + {1
\over n^2} {q_{n} / q'_{n}}$. Hence every term in the sum above is
constant on  $D_n^{iq_{n}
q'_{{n}-1}}, $ equal to $\epsilon_n$ if $l \leq i$ and to $0$ if $l > i$ and the proposition follows.    \carre

\vspace{0.2cm}

Since for every $y$ at most one of
 the functions $ \theta \left(n^2 {q'_n / q_n} (y - i {q_n / q'_n}) \right) $ is not locally constant we obtain the following straightforward
estimates for $\kappa_n$:
\begin{lemm} \label{kappacontrol} For any $p \geq 1$, we have 
\begin{eqnarray} \label{kk11}
\| \kappa_n \| &\leq& r_n \epsilon_n \\
{\| {\kappa}_{n} \|}_{C^p} &\leq&   2{\|\theta\|}_{C^p} n^{2p} { \left( {q'_{n} \over q_{n}} \right)}^p \epsilon_n.  
\end{eqnarray}
\end{lemm}

To complete our construction of the function $X_n$, we still have  to adjust $\kappa_n$ in
order to have a function that satisfies (X.3) not only on the sets $B_n^{i q_{n}
q'_{n-1}} \subset D_{n}^{ i q_{n}q'_{n-1} }  \subset R_n^0$, but on all the levels ${B}_n^h$ of the tower $\GG (B_n^0, h_n) $
(see the definitions in Section \ref{vvvspecialtower}).  Define on the real line the functions 
 \begin{eqnarray} 
\label{biu} \nu_n (x)  &=&  \theta \left( n^2 q_{n}x \right) 
- \theta \left( n^2 q_{n} ( x-  {1 \over q_{n}} + { 2 \over n^2 q_{n}} )  \right), \\ \upsilon_n(y)   &=& \theta \left( n^2 q'_{n-1} y + 2 \right) 
 - \theta \left( n^2 q'_{n-1} (y-  {1 \over q'_{n-1}} + {3 \over n^2 q'_{n-1}})  \right).
 \label{bip} \end{eqnarray}
We then have
\begin{eqnarray*}
\nu_n(x) &=& 0 \ \ {\rm for \ \ } x \leq 0, {\rm \ and \ } x \geq { 1 \over q_{n}} - {1 \over n^2 q_{n}}, \\
\nu_n(x) &=& 1 \ \ {\rm for \ \ }   { 1 \over n^2 q_{n}} \leq x \leq { 1 \over q_{n}} - {2 \over n^2 q_{n}},
\end{eqnarray*}
and
\begin{eqnarray*}
\upsilon_n(y) &=& 0  \ \ {\rm for \ \ } y \leq - {2 \over n^2 q'_{n-1}}, {\rm \ and \ } y \geq { 1 \over q'_{n-1}} - {2 \over n^2 q'_{n-1}}, \\
\upsilon_n(y) &=& 1  \ \ {\rm for \ \ } - {1 \over n^2 q'_{n-1}} \leq y \leq { 1 \over q'_{n-1}} - {3 \over n^2 q'_{n-1}}.
\end{eqnarray*}

Hence we can consider the restriction of the function $\nu_n(x) \upsilon_n(y)$ on the rectangle $[0, 1 /q_n] \times [-2 / n^2 q'_{n-1}, 1 / q'_{n-1} - 2 / n^2 q'_{n-1}]$ and extend it by $0$ to a $C^\infty$ function on the two torus. With $\overline{R}_n^0$ as in Definition \ref{RR} we have that the latter function is constant equal to 1  on $\overline{R}_n^0$ and to 0 on all the other $\overline{R}_n^j$, $1 \leq j \leq q_{n}q'_{n-1} -1$. It is also  easy to see that the functions 
$ \nu_n ( x - l {p_{n} / q_{n}}) \upsilon_n ( y - l {p'_{{n}-1}
/ q'_{{n}-1}}) $ are  equal to 1 on $\overline{R}_n^l$ and to 0 on all the other $\overline{R}_n^j$, $0 \leq j \leq q_{n}q'_{n-1} -1$.

 Having Corollary \ref{vvvescalier} in mind and the definition of $\beta_{n,j}$ in Lemma \ref{vvvzabb}, we set
  \begin{eqnarray} \label{1xn} \tilde{X}_n (x,y) :=  \sum_{j=0}^{q_{n} q'_{{n} - 1} -1 } \kappa_n ( y -
{\beta}_{n,j})  \nu_n ( x - j {p_{n}
\over q_{n}}) \upsilon_n ( y - j {p'_{{n}-1}
\over q'_{{n}-1}}) \end{eqnarray}
where $\kappa_n$, $\nu_n $ and $\upsilon_n$ are defined in (\ref{bas}), (\ref{biu}) and (\ref{bip}).
The latter function is of class $C^\infty$ on  ${\TT}^2$ and is equal to $\kappa_n ( y - \beta_{n,j} ) $ on $\overline{R}_n^j$. Hence, Corollary \ref{vvvescalier} and Lemma \ref{kkkkk} imply that
\begin{prop} \label{vvvrangunadroite} With $h_n$ and $B_n^p$ as in Definition \ref{B}, the function  $\tilde{X}_n$ is constant on every ${B}_n^p$,
for all  $p 
\leq h_n$. More precisely, for every $0 \leq j \leq q_{n} q'_{n-1} -1$, for every $i \leq (1 - {2 \over n}) r_n$, we have 
\begin{eqnarray} \label{tutto}{{\tilde{X}}_n} { |}_{ B_n^{j+i q_{n} q'_{n-1}}} \equiv  i \epsilon_n.\end{eqnarray}
\end{prop}

The following estimates follow immediately from (\ref{kk11})--(\ref{1xn})

\begin{prop} \label{vvvchicontrol}  We have for $n$ large enough 
\begin{eqnarray} 
\|\tilde{X}_n \| \leq r_n \epsilon_n, \end{eqnarray}
and for $(r,p) \neq (0,0)$ we have for some constant $c(r,p)$ 
\begin{eqnarray} \|D^r_x D^p_y \tilde{X}_n \| \leq c(r,p) n^{2r+2p} q_n^r {\left({q'_n \over q_n} \right)}^p \epsilon_n.
\end{eqnarray}
 
\end{prop}

\vspace{0.2cm} 

\subsection{The function ${X_n}$.}

$\\$

Propositions \ref{vvvrangunadroite} and \ref{vvvchicontrol} enclose  the properties (X.1), (X.2) and (X.3)
required for $X_n$. Nevertheless the function $\tilde{X}_n$ may fail to satisfy the uniform bounds required in (X.4) and (X.5). The presence of high frequencies in the Fourier
expansion of $\tilde{X}_n$  may indeed put in jeopardy such control. To overcome this problem we turn
to truncations. The function $\tilde{X}_n$ defined in (\ref{1xn}) being of class $C^\infty$ we consider its Fourier coefficients $X_{l,j}^n$ and define

 \begin{eqnarray} \label{1xntronque} {X}_n (x,y) :=  
\sum_{l^2 + j^2 \leq {q'_n}^4; (l,j) \notin (\ZZ q_n, \ZZ^* q'_n)} X_{l,j}^n \chi_{l,j}, 
\end{eqnarray}
where for $(l, j) \in \ZZ^2$ 
$$\chi_{l,j}(x,y) := e^{i2 \pi (lx+jy)}.$$
It is for the function $X_n$ that we want to check the properties of Proposition \ref{X}.
\begin{lemm} \label{1xnn} For any $r \geq 1$, we have for $n$ sufficiently large
\begin{eqnarray}  \| {X}_n  - \tilde{X_n} \| \leq {1 \over q_nq'_n}, \label{1C2} \\  
 {\| D_y({X}_n  - \tilde{X}_n)  \|} &\leq& {1 \over q_n^2}, \label{1C1} \\
  \|D_x^r  ({X}_n  - \tilde{X}_n) \|  &\leq&  q_n^{r-1} \epsilon_n. \label{1C0}  
\end{eqnarray}
\end{lemm}

\vspace{0.15cm}

\noindent {\sl Proof.} By definition 
$$\tilde{X}_n - X_n =   \sum_{l^2 + j^2 > {q'_n}^4} X_{l,j}^n \chi_{l,j} + 
\mathop{\sum_{l^2 + j^2 \leq {q'_n}^4}} \limits_{(l,j) \in (\ZZ^*q_n, \ZZ^*q'_n)} X_{l,j}^n \chi_{l,j} +  \sum_{0<|p| \leq q'_n} X_{0,pq'_n}^n \chi_{0,pq'_n}, $$ 
we will denote in this proof by A, B, C the first, second and third term in the latter expression of $\tilde{X}_n - X_n$.

Since ${(l^2+ j^2)}^2  |X^n_{l,j}| \leq {\| \tilde{X}_n\|}_{C^4} $  and ${(pq'_n)}^2|X_{0,pq'_n}| \leq {\| \tilde{X}_n\|}_{C^2} $ we have 
\begin{eqnarray*}  \| A+ C \| &\leq&    {\| \tilde{X}_n\|}_{C^4}  \sum_{l^2 + j^2 > {q'_n}^4} { 1 \over {(l^2+ j^2)}^2 } + { { \| \tilde{X}_n\|}_{C^2} \over {q'_n}^2 }  \sum_{0<|p| \leq q'_n} {1 \over p^2} \\
&\leq& c_1    {{\| \tilde{X}_n\|}_{C^4} \over {q'_n}^4}+ c_2 {{\| \tilde{X}_n\|}_{C^2} \over {q'_n}^2} \\
&=& o({1 \over q_nq'_n})
\end{eqnarray*}
from Proposition \ref{vvvchicontrol} and our choice of $\epsilon_n$ in (\ref{epsn}). Likewise
\begin{eqnarray*}  \| B \| &\leq&    \sum_{(l,j) \in (\ZZ^* q_n, \ZZ^* q'_n)}  { \| D_x^2 D_y^2 \tilde{X}_n\| \over l^2j^2 }  \\
&\leq& {C \over q_n^2 {q'_n}^2}  {\| D_x^2 D_y^2 \tilde{X}_n\|} \\
 &=& o({1 \over q_nq'_n}) \end{eqnarray*} 
again from Proposition \ref{vvvchicontrol} and our choice of $\epsilon_n$ in (\ref{epsn}). We have proved (\ref{1C2}).

In the same way as above we have 
\begin{eqnarray*}  \|D_y(\tilde{X}_n - X_n)\| 
&\leq& c_3    {{\| \tilde{X}_n\|}_{C^5} \over {q'_n}^4}+ c_4 { \| D_x^2 D_y^3 \tilde{X}_n\| \over q_n^2 {q'_n}^2}  + c_5 { { \| \tilde{X}_n \|}_{C^3} \over {q'_n}^2 } \\
&=& o({1 \over q_n^2}).
\end{eqnarray*}
 
Finally, for any $r \geq 1$ we have 
\begin{eqnarray*}  \|D_x^r(A)\|
&\leq& c_5    { {\| D_x^r \tilde{X}_n\|}_{C^4} \over {q'_n}^4 } \\
&=& o(q_n^{r-1} \epsilon_n)
\end{eqnarray*}
from Proposition \ref{vvvchicontrol} and our choice of $\epsilon_n$ in (\ref{epsn}), while we also have
\begin{eqnarray*}  \|D_x^r(B)\|
&\leq& c_6    { {\| D_x^{r+2} D_y^2 \tilde{X}_n\|} \over q_n^2 {q'_n}^2}  \\ 
&=& o(q_n^{r-1} \epsilon_n). \end{eqnarray*}  \carre

\subsection{Checking the properties of Proposition  \ref{X} for $X_n$.} In light  of Lemma \ref{1xnn}, the Propositions \ref{vvvrangunadroite} and \ref{vvvchicontrol} yield (X.3), (X.1) and (X.2). It remains to give the

\vspace{0.2cm}

\noindent {\sl Proof of (X.4) and (X.5).}   We will need the following Lemma
 
\begin{lemm} \label{ibound} We have for any $m \in \NN$
\begin{eqnarray} \label{ibound1} \|S_m \chi_{l,j}\| \leq \inf \left( m , {1 \over 2 ||| l\a+j\a'||| } \right) \end{eqnarray}
where $|||.|||$ denotes the  distance to the closest integer. 
\end{lemm}

\noindent {\sl Proof.} We have
$$S_m \chi_{l,j} (x,y) = U_{l,j,m} \chi_{l,j}(x,y),$$
where 
 $$ U_{l,j,m} = {1 -  e^{i 2 \pi m (l \a + j \a')} \over  1 -
e^{i 2 \pi  (l \a + j \a')}}.$$
Clearly 
\begin{eqnarray*}
|U_{l,j,m}| &=& \left| {\sin (\pi m(l\a + j \a')) \over \sin (\pi (l\a+j\a'))} \right| \\ 
&\leq& {1 \over | \sin (\pi(l\a+j\a')) |},
\end{eqnarray*}
but $|\sin(\pi(l\a+j\a')| = |\sin(\pi(\tn l\a+j\a' \tn)| \geq 2 \tn l\a+j\a' \tn.$ \carre
 
\vspace{0.2cm}

We have
\begin{eqnarray*} S_m X_n (x,y) = \sum_{l^2+j^2 \leq {q'_n}^4; (l,j) \notin (\ZZ q_n, \ZZ^* q'_n)} X^n_{l,j}   S_m \chi_{l,j} \\
= \sum_{l^2+j^2 \leq {q'_n}^4; (l,j)  \notin (\ZZ q_n, \ZZ q'_n)} X^n_{l,j}   S_m \chi_{l,j} + \sum_{|l| \leq {{q'_n}^2 \over q_n}}  S_m \chi_{lq_n,0} \end{eqnarray*}

From (\ref{reduites2}) in Section \ref{notations} and Properties (\ref{10000})--(\ref{30000}) in our choice of $\a,\a'$ in Section \ref{choice} we deduce that 

\begin{itemize} 

\item For $ |l| \leq {q_n'}^2$ and $|j| \leq {q'_n}^2 $ such that $|j| \notin \NN q'_n$ or $|l| \notin \NN q_n$ we have  
\begin{eqnarray} \label{crucial} \tn l\a+j\a' \tn \geq {1 \over 2 q_{n} q'_n}. \end{eqnarray} 
 
\item For $ |l| < {q_{n+1}}$ 
\begin{eqnarray} \label{cc} \tn l \a \tn \geq {1 \over 2 q_{n+1} }. \end{eqnarray}
 
\end{itemize}

Using Lemma \ref{ibound} we hence get for $r \geq 1$
\begin{eqnarray*} 
\| D_x^r S_m X_n (x,y)\| 
&\leq& q_n q'_n \sum_{l^2+j^2 \leq {q'_n}^4}   
{(2 \pi|l|)}^r |X^n_{l,j}|   \\
 \ & \ &  \   + q_{n+1} \sum_{0< |l| \leq {q'_n}^2}  {(2 \pi|l|)}^r |X^n_{l,0}|  \\ 
 &\leq& c q_n {q'_n}^5  \|D_x^r \tilde{X}_n \| + q_{n+1} \|D_x^{r+2} \tilde{X}_n \| \sum_{0< |l| \leq {q'_n}^2 } {1 \over l^2},  
\end{eqnarray*}
which yields (X.4) of Proposition \ref{X} due to Proposition \ref{vvvchicontrol}.

Similarly, to get (X.5) we write for $p \geq 1$ 
\begin{eqnarray*} 
\| D_y^p S_m X_n (x,y)\| 
&\leq& q_nq'_n \sum_{l^2+j^2 \leq {q'_n}^4}  {(2 \pi|j|)}^p |X^n_{l,j}|  \\
&\leq&  q_n  {q'_n}^5 \|D_y^p \tilde{X}_n\| 
\end{eqnarray*}
which implies (X.5) due to Proposition \ref{vvvchicontrol}.       \carre

\section{\sc Construction of $Y_n$.} \label{psiconstruct}

\subsection{Construction of a first function $\tilde{Y}_n$.}

Define, for  $y \in [0,{1\over q'_{{n}-1}}]$ the following function 
\begin{eqnarray*}  \tilde{\phi}_n(y) :=  \theta \left( n^3 q'_{{n}-1} y - n^3 +n^2  \right) - \theta \left( n^3 q'_{{n}-1}  y - n^3 +n  \right). \hspace{0.2cm} \end{eqnarray*} 
Since $\theta$ is increasing and $\theta(s) = 0$ for $s \leq 0$, and $\theta(s) = 1$ for
$s \geq 1$, it is easy to check that 
\begin{eqnarray} \label{vvvphileq1}
0 \leq \tilde{\phi}_n(y) \leq 1, \end{eqnarray} 
and that
\begin{eqnarray} \label{vvvphi0} 
 \tilde{\phi}_n(y) = 0, { \rm \ \ if \ \ }  y \in \left[0,(1 - {1 \over
n}){1\over q'_{{n}-1}} \right] \cup  \left[(1 - {1 \over
n^3}){1\over q'_{{n}-1}}, {1\over q'_{{n}-1}} \right] ,    \end{eqnarray}
while
\begin{eqnarray}
\tilde{\phi}_n(y) = 1, { \rm \ \  if \ \   }  y \in \left[(1 - {1 \over
n}+ {1\over n^2}){1\over q'_{{n}-1}}, (1 - {1 \over
n^2}){1\over q'_{{n}-1}}      \right]. \label{vvvphi1} \end{eqnarray}

Due to  (\ref{vvvphi0}) it is possible to extend  $\tilde{\phi}_n$ to the
circle as a  $C^{\infty}$  periodic function with period ${1 / q'_{{n}-1}}$ with the following  estimate for any $p \in \NN$ 
\begin{eqnarray} \label{propo}
 \| D^p_y \tilde{\phi}_n  \| =  {\|\theta\|}_{C^p}  n^{3p} {(q'_{n-1})}^p. \end{eqnarray}

 \vspace{0.2cm}

Define now on the two torus the function 
$$\tilde{Y}_n(x,y) := { - \cos (  2 \pi q_{n}x) \over e^{q_{n} } } \tilde{\phi}_n(y).$$

The function $\tilde{Y}_n$ is of class $C^{\infty}$, and satisfies 
\begin{prop} \label{vvvprepsicontrol} For $r, p \in \NN$, we have
\begin{eqnarray*}
\| D^r_x D^p_y \tilde{Y}_n  \| =   {\|\theta\|}_{C^p} {(2 \pi)}^r n^{3p} {(q'_{n-1})}^p {q_n^r \over e^{q_{n}}}. 
\end{eqnarray*}
\end{prop}

In preparation for (Y.3) and (Y.3') we have

\begin{prop} \label{vvvpsistretch}
For $x$ such that $\lbrace q_{n} x \rbrace  \in [{1 \over n}, {1 \over 2} -{1 \over n}] \cup  [{1 \over 2} + {1 \over n}, 1 -{1 \over n}]$, and  $m \in [q'_{n}, 2 q_{n+1} / {(n+1)}^2 ]$, we have
\begin{eqnarray*}
\left|   D_x S_m \tilde{Y}_n  (x,y) \right| &\geq&
  {\pi \over n^2} {q_n \over  e^{q_{n}} } m,     \\  
\left\| D^2_x   S_m \tilde{Y}_n 
\right\| &\leq& 4 \pi^2 { q_n^2 \over
  e^{q_{n}} } m. \end{eqnarray*}
Given $\eta>0$, the same inequalities above hold for sufficiently large $n$ for $m \in [q'_n / 2n^2, q'_n]$ if we restrict $y$ to  $1-m/q'_n+\eta \leq \lbrace q'_{n-1} y \rbrace \leq 1 - \eta$.
\end{prop}$ $

\noindent {\sl Proof. } We have
$${D_x \tilde{Y}_n  } (x,y) = 2 \pi q_n {\sin( 2 \pi
q_{{n}}x) \over e^{q_{{n}}}} \tilde{\phi}_n(y).$$
 Assume that $ \lbrace q_{n} x \rbrace \in [ {1\over n }, { 1 \over 2 } -
{1\over n }]$, the other case  being similar. For $k \leq m \leq 2 q_{n+1} / {(n+1)}^2$ (\ref{reduites1}) implies that 
 $$ \lbrace q_{n}(x + k \a) \rbrace  \in   [{1 \over 2 n }, { 1 \over 2 }
- {1 \over 2 n }],$$ 
hence, $ \sin ( 2 \pi q_{{n}} ( x + k \a ) ) \geq \sin({\pi / n}) \geq {2 / n} $ and because $\tilde{\phi}_n $ is positive this implies
$$ { D_x \tilde{Y}_n }  (x+k\a,y+ k\a')   \geq \tilde{\phi}_n(y +
  k\a')
  {4 \pi q_{n} \over n e^{q_{n}} }
  \geq 0. $$
In light of (\ref{vvvphi1}), we will finish if we prove that for every $y$ and for every $m \geq q'_{n}$ there is more than $ m  / 4n $ integers $k \leq m$ such that 
$ \lbrace q'_{{n}-1}( y + k \a') \rbrace  \in  \left[1 - {1 \over
n}+ {1\over n^2}, 1 - {1 \over
n^2}      \right] $. The latter follows if  for every $y$, there is at least ${q'_n \over 2n}$ integers $k \leq q'_n$ satisfying the desired property. This in  turn follows from the good approximation of  $R_{\a'}$ by $R_{p'_n / q'_n}$.  The proof in the case $m \in [q'_n /2n^2,q'_n]$ follows in the same way.

\noindent To obtain the inequality involving the second derivative we just bound the
cosine by $1$ and use $\| \tilde{\phi}_n \| \leq 1 $. \carre 

\vspace{0.2cm}

In preparation for (Y.2) we have
\begin{prop}\label{vvvprepsirankone}  With $h_n$ and $B_n^0$ as in Definition \ref{B}    we have  for any $m \leq h_n$ that $  S_m \tilde{Y}^{n}$
is identically zero on ${B}_n^0$. 
\end{prop}

\noindent {\sl Proof. } Given $h \leq m \leq h_n$, let $i \leq (1 - 2 /n)r_n$ and $0 \leq j \leq q_nq'_{n-1}-1$ be such that $h= j+ iq_nq'_{n-1}$. From (\ref{hh2}) in Section \ref{vvvspecialtower} we have that for $(x,y) \in R_{{\a,\a'}}^h ({B}_n^0) $ 
\begin{eqnarray*} y &\in& \left[ j {p'_{n-1} \over q'_{n-1}}, j {p'_{n-1} \over q'_{n-1}} + (i+2){q_n \over q'_n} \right] \\
&\subset& \left[   j {p'_{n-1} \over q'_{n-1}}, j {p'_{n-1} \over q'_{n-1}} + (1 - {1 \over n}) {1 \over q'_{n-1}} \right], \end{eqnarray*}
hence  $\tilde{\phi}_n(x,y)= 0$ from (\ref{vvvphi0}). \carre

\subsection{The function  ${Y}_n$.}

 \vspace{0.2cm}

The function $\tilde{Y}_n$ has all the required properties by Proposition \ref{Y} except for the uniform bounds on the Birkhoff sums of the derivatives.
As we did in the last section we replace $\tilde{\phi}_n$ by 
$${\phi}_n (y) = \sum_{ j= - q'_n +1}^{ q'_n -1 } {\phi}^n_j e^{i 2 \pi j y},$$ 
where the  ${\phi}^n_j$ are the Fourier coefficient of $\tilde{\phi}_n$. We then let  
$${Y}_n(x,y) :=  - { \cos ( 2 \pi q_nx ) \over e^{q_n}} {\phi}_n(y).$$
The truncation here is less delicate than in the definition of $X_n$ due to the fact that the sequence $Y_n$ converges to zero in the $C^\infty$ norm.

\subsection{Checking the properties of Proposition \ref{Y} for $Y_n$.} 

From (\ref{propo})  it is easy to see that
for any $p \in \NN$, we have for $n$ large enough
\begin{eqnarray*} {\| \tilde{\phi}_n  - {\phi}_n  \|}_{C^p} &\leq& \sum_{|j| \geq q'_n} {(2 \pi |j|)}^p |\phi_j^n| \\
&\leq&  {\| \tilde{\phi} \|}_{C^{p+4}} \sum_{|j| \geq q'_n } {1 \over j^4} \\
&=&  
 o({1 \over
 {q'_n}^2 }) \end{eqnarray*}
Combined with Propositions \ref{vvvprepsicontrol}, \ref{vvvpsistretch} and \ref{vvvprepsirankone} the above yields (Y.1), (Y.2), (Y.3) and (Y.3')  for $Y_n$. It remains to give the

\noindent {\sl Proof of (Y.4) and (Y.5).}
The proof  is  similar yet easier than that for $X_n$:  For any $m \in \NN$ we have 
$$ S_m Y_n(x,y) =  \sum_{|j| <q'_n}  {{\phi}^n_j \over 2 e^{q_n}} S_m (\chi_{q_n,j} + \chi_{-q_n,j}).  $$
 Since $\tn q_n \a \tn < 1 /q_{n+1}$ while $\tn j \a' \tn > 1 / 2 q'_n$ for any $0< |j| < q'_n$ then for such $j$  we have $\tn  \pm q_n\a+j\a' \tn \geq 1 / 4 q'_n$, hence Lemma \ref{ibound} implies 
\begin{eqnarray*}  \|D_y^p S_m Y_n (x,y)\| &\leq& 2 {q'_n \over e^{q_n}}  \sum_{0< |j| <q'_n}  {(2 \pi |j|)}^p {|{\phi}^n_j|  }   \\
&\leq& C {q'_n \over e^{q_n}}  {\|\tilde{\phi}_n\|}_{C^{p+2}} \end{eqnarray*}  
which yields (Y.5) due to (\ref{propo}).

For $r \in \NN$ we have 
 \begin{eqnarray*}  D_x^r S_m Y_n (x,y) &=& 
 \sum_{0< |j| <q'_n}  {{\phi}^n_j  \over 2 e^{q_n}} {(i 2 \pi q_n)}^r S_m (\chi_{q_n,j} + \chi_{-q_n,j}) \\
&+&   {{\phi}^n_0 \over 2 e^{q_n}} {(i 2 \pi q_n)}^r S_m (\chi_{q_n,0} + \chi_{-q_n,0}). \end{eqnarray*} 
As in our proof of (Y.5) the first term is uniformly bounded away from  $q'_n$; while the second term is bounded by $q_{n+1} {(2 \pi q_n)}^r / e^{q_n}$ since $0 \leq \phi^n_0 \leq \|\tilde{\phi}_n\| \leq 1$ and $\|S_m \chi_{ \pm q_n,0}\| \leq 1 / (2 \tn q_n \a \tn) \leq q_{n+1}$. Hence (Y.4) is proved. \carre

\vspace{1cm} 

\noindent {\sl Acknowledgments.} I wish to express my gratitude to the late Michael Herman for encouraging me to do this work.  I also wish
to thank  Jean-Paul Thouvenot for fruitful conversations since the early stage of the paper and Jean-Christophe Yoccoz for suggesting simplifications in the construction.  

\bibliographystyle{plain}                  
\bibliography{ii} $ $

\end{document}